\documentclass[12pt]{amsart}

\usepackage[utf8]{inputenc}
\usepackage[T1]{fontenc}
\usepackage[english]{babel}

\usepackage{amsfonts}
\usepackage{amsmath}
\usepackage{amssymb}
\usepackage{amstext}
\usepackage{amsthm}

\usepackage{geometry}
\usepackage{fullpage}

\usepackage[shortlabels]{enumitem}
\usepackage{graphicx}
\usepackage[hidelinks]{hyperref}
\usepackage{mathrsfs}
\usepackage{mathtools}

\usepackage{tikz}
\usetikzlibrary{calc,decorations.pathreplacing}

\allowdisplaybreaks

\newcommand{\cA}{\mathcal{A}}
\newcommand{\Ex}{\mathbb{E}}
\DeclareMathOperator{\Var}{Var}
\DeclareMathOperator{\SYT}{SYT}

\newtheorem{theorem}{Theorem}[section]
\newtheorem{lemma}[theorem]{Lemma}
\newtheorem{proposition}[theorem]{Proposition}
\newtheorem{corollary}[theorem]{Corollary}

\theoremstyle{definition}

\theoremstyle{remark}
\newtheorem{remark}[theorem]{Remark}

\theoremstyle{plain}

\author{Ran Pan}
\address{Department of Mathematics, University of California San Diego, La Jolla, CA 92093-0112, USA}
\email{ran.pan.math@gmail.com}
\thanks{}

\author{Dun Qiu}
\address{Center for Combinatorics, LPMC, Nankai University, Tianjin 300071, P. R. China}
\email{qiudun@nankai.edu.cn}
\thanks{Dun Qiu is supported in part by the National Natural Science Foundation of China (No. 12671390),
the Natural Science Foundation of Tianjin (No. 24JCZDJC01390), and the Fundamental Research Funds for the Central Universities (No. 105-63263093).}

\author{Weichen Wang}
\address{Center for Combinatorics, LPMC, Nankai University, Tianjin 300071, P. R. China}
\email{wangweichenc@163.com}
\thanks{}

\author{Xinhui Wang}
\address{Center for Combinatorics, LPMC, Nankai University, Tianjin 300071, P. R. China}
\email{wangxinhui2018@163.com}
\thanks{}

\title{Entry and position statistics in alternating permutations and Young tableaux}

\date{\today}



\begin{document}

\begin{abstract}
We study entry and position statistics in alternating permutations and standard Young tableaux using linear extensions and generating functions. For alternating permutations, we give explicit formulas for the mean and variance of the first entry and of the position of any fixed value. For standard Young tableaux, we express the sum of the lower-left entries as the number of tableaux of an associated skew shape. For two-row shapes, we determine the distribution of the entry in each cell of the second row, and for the shape $(n,n)$ we give a closed formula for its sum over all tableaux, which recovers a limit formula of Stanley. Through a bijection of Lewis, the results for the shape $(n,n)$ describe the entries of a random $321$-avoiding down-up permutation of even length.
\end{abstract}

\maketitle

\section{Introduction}

Alternating permutations are a classical family in enumerative combinatorics. Their enumeration, initiated by Andr\'e~\cite{andre1881permutations}, gives the Euler numbers. Entringer~\cite{entringer1966combinatorial} refined this enumeration by recording the first entry. This refinement provides a natural starting point for studying the distribution of entries in random alternating permutations. Further aspects of alternating permutations and their connections with other combinatorial structures are discussed in Stanley's survey~\cite{stanley2009survey}. Entringer numbers have been studied from several points of view; see Gelineau, Shin and Zeng~\cite{gelineau2011bijections} for bijections between the families counted by them, and Foata and Han~\cite{foata2014seidel} for the bi-Entringer numbers, which count alternating permutations by both the first and the last entry. The limiting distribution of the first entry of a random alternating permutation was found by Simons and Yao~\cite{simons1996sine}. The means and variances of other permutation statistics have been studied extensively; Knuth~\cite{donald1999art} treats inversions and runs, and B\'ona~\cite{bona2007copies,bona2010three} treats the number of copies of a pattern and three notions of monotone subsequences.

In this paper, we study entry and position statistics through linear extensions of finite posets. An alternating permutation can be represented by an order-preserving labeling of a zigzag poset, that is, by a linear extension, and a standard Young tableau is a linear extension of the poset of cells of its shape. The first entry of an alternating permutation and the entry in a prescribed cell of a tableau are therefore both labels of fixed poset elements in a uniformly random linear extension; the average of such a label, the average height of the element, was studied by Winkler~\cite{winkler1982average}. Adjoining elements to the poset converts weighted sums of such labels into ordinary counts of linear extensions. We use this construction to determine the expectation and variance of the first entry, and the sum of the lower-left entries of the standard Young tableaux of a given shape. We also use a bivariate exponential generating function and an argument based on interchanging consecutive values to study the position $X_j(\sigma)=\sigma^{-1}(j)$ of a prescribed value $j$.

For a uniformly random up-down alternating permutation $\sigma$ of length $n$, we show that
\[
\Ex[\sigma_1]
=
n+1-\frac{E_{n+1}}{E_n},
\]
where $E_n$ denotes the $n$-th Euler number, together with a formula for the variance (Theorem~\ref{thm:first-entry}). For permutations of even length $2n$, we show that for every $1\le j\le2n$,
\[
\Ex[X_j]
=
n+\frac12
-
\frac{2n-2j+1}{2}\,\frac{E_{2n-2}}{E_{2n}}
\]
and
\[
\Var(X_j)
=
\frac{4n^2-1}{12}
+
\frac{-2n^2+(6j-3)n-3j(j-1)-1}{3}\,
\frac{E_{2n-2}}{E_{2n}}
-
\frac{(2n-2j+1)^2}{4}
\left(\frac{E_{2n-2}}{E_{2n}}\right)^2
\]
(Theorem~\ref{thm:general-position}(a)). For odd length $2n+1$, every value has expected position $n+1$, and
\[
\Var(X_j)=\frac{n(n+1)}{3}+\left(\frac{2n(2n-1)}{3}-(j-1)(j-2)\right)\frac{E_{2n-1}}{E_{2n+1}}
\]
(Theorem~\ref{thm:general-position}(b)). In both cases, for fixed $j$ and length $N$, the mean and the variance of $X_j$ are $(N+1)/2+O(N^{-1})$ and $(N^2-1)/12+O(1)$, the mean and the variance of a uniformly random position.

Entries of random standard Young tableaux have been studied from several directions. McKay, Morse and Wilf~\cite{mckay2002distributions} determined the limiting probability that prescribed entries occupy prescribed cells in a random tableau with $n$ cells, with no restriction on the shape. For a fixed shape, Riegler and Neumann~\cite{riegler2016playing} gave a product formula for the mean of the entry in the cell $(2,1)$, and Ekhad and Zeilberger~\cite{ekhad2024experimenting} computed the distribution of the entry in a given cell by symbolic computation. Let $\SYT(\lambda)$ denote the set of standard Young tableaux of shape $\lambda$, let $f(\lambda)$ denote the number of standard Young tableaux of shape $\lambda$, and let $T(i,j)$ denote the entry in row $i$ and column $j$ of a tableau $T$. For a shape $\lambda$ with $k$ parts, adding one cell to the left of the lower-left cell shows that the sum of the entries $T(k,1)$ over $\SYT(\lambda)$ is the number of standard Young tableaux of a skew shape (Theorem~\ref{thm:general-shape}). For the rectangular shape $n^k=(n,\ldots,n)$ with $k$ parts, this gives
\[
\sum_{T\in\SYT(n^k)}T(k,1)
=
f(n+1,n^{k-1}),
\]
and the mean of $T(k,1)$ is $k(nk+1)/(n+k)$ (Corollary~\ref{cor:rect}). For two-row shapes $(m,n)$ with $m\ge n\ge1$, we determine the distribution of the entry in each cell of the second row. For a uniformly random $T\in\SYT(m,n)$ and $1\le k\le n$, we obtain the following explicit formula for the probability that $T(2,k)=a$:
\[
\Pr(T(2,k)=a)
=
\frac{
	\left[\binom{a-1}{k-1}-\binom{a-1}{k-2}\right]
	\left[\binom{m+n-a}{n-k}-\binom{m+n-a}{n-a+k-1}\right]
}{
	\binom{m+n}{n}-\binom{m+n}{n-1}
},
\qquad 2k\le a\le m+k,
\]
and $\Pr(T(2,k)=a)=0$ for all other integers $a$ (Theorem~\ref{thm:two-row-dist}). For the shape $(n,n)$, we also show that
\[
\sum_{T\in\SYT(n,n)}T(2,k)
=kC_{n+1}-\sum_{i=1}^{k-1}(k-i)\,C_{i-1}C_{n+1-i},
\]
where $C_n$ denotes the $n$-th Catalan number (Theorem~\ref{thm:second-row-sums}); dividing by $C_n$ and letting $n\to\infty$ recovers a limit formula of Stanley for the mean of $T(2,k)$.

Mansour~\cite{mansour2003restricted} initiated the study of pattern avoidance in alternating permutations, and Lewis~\cite{lewis2011pattern} related it to standard Young tableaux. For permutations avoiding a pattern of length three, Miner and Pak~\cite{miner2014shape} studied the asymptotic distribution of the positions of the values, and the number of $123$- and $132$-avoiding permutations with a given first entry is a ballot number; see Miner and Pak~\cite{miner2014shape} and E\u{g}ecio\u{g}lu, Gaiser and Yin~\cite{egecioglu2023leading}. Reading the columns of a standard Young tableau of shape $(n,n)$ gives a bijection onto the $321$-avoiding down-up permutations of length $2n$ (Lemma~\ref{lem:reading-word}), which is a form of a bijection of Lewis. Through it, our results for the shape $(n,n)$ give the mean of every entry of a random $321$-avoiding down-up permutation $w$ of length $2n$, and the distribution of every entry in an odd position (Theorem~\ref{thm:pattern}). In particular, $\Ex[w_1]=C_{n+1}/C_n$, whereas $\Ex[w_1]=E_{2n+1}/E_{2n}$ for a random down-up permutation of length $2n$ (Corollary~\ref{cor:down-up-first}).

The paper is organized as follows. Section~2 introduces the notation and recalls some preliminary results. Section~3 derives formulas for the expectation and variance of the first entry in alternating permutations. Section~4 studies the position statistics of fixed values in alternating permutations. Section~5 presents results on entry statistics in standard Young tableaux, with particular attention to rectangular and two-row shapes, and applies the results for the shape $(n,n)$ to pattern-avoiding alternating permutations.

\section{Background}

Let $S_n$ denote the symmetric group on $n$ letters. In this paper, we adopt one-line notation $\sigma=\sigma_1\sigma_2\cdots\sigma_n$ for permutations. A permutation $\sigma\in S_n$ is \emph{up-down alternating} if $\sigma_1<\sigma_2>\sigma_3<\sigma_4>\cdots$, and \emph{down-up alternating} if $\sigma_1>\sigma_2<\sigma_3>\sigma_4<\cdots$. For example, $2746153$ is up-down alternating, while $4152736$ is down-up alternating. We denote the set of all up-down (down-up) alternating permutations in $S_n$ by $\cA_n$ ($\cA'_n$). The \emph{complement} map $\sigma\mapsto\sigma'$, where $\sigma'_i=n+1-\sigma_i$, is a bijection from $\cA_n$ to $\cA'_n$, so both sets have the same cardinality $E_n=|\cA_n|=|\cA'_n|$, the $n$-th \emph{Euler number}. The first few Euler numbers are $1,1,1,2,5,16,61,\ldots$; see sequence A000111 in the OEIS~\cite{OEIS}. It is well known that the exponential generating function of the Euler numbers is
\begin{equation}\label{eq:andre}
\sum_{n\geq 0}E_n\frac{x^n}{n!}=\sec x+\tan x;
\end{equation}
this goes back to Andr\'e~\cite{andre1881permutations}. Moreover, separating the even and odd parts of \eqref{eq:andre} gives
\begin{equation}\label{eq:sectan}
\sum_{n\geq 0}E_{2n}\frac{x^{2n}}{(2n)!}=\sec x,\qquad
\sum_{n\geq 0}E_{2n+1}\frac{x^{2n+1}}{(2n+1)!}=\tan x.
\end{equation}
The Euler numbers have the asymptotic behavior (see Stanley~\cite{stanley2009survey})
\begin{equation}\label{eq:asym}
\frac{E_n}{n!}
=
\frac{4}{\pi}\left(\frac{2}{\pi}\right)^{n}
+O\left(\left(\frac{2}{3\pi}\right)^{n}\right).
\end{equation}
In particular,
\begin{equation}\label{eq:ratio}
\frac{E_{n+1}}{E_n}
=
\frac{2}{\pi}(n+1)\bigl(1+O(3^{-n})\bigr),
\qquad\text{and}\qquad
\lim_{n\to\infty}\frac{E_n}{nE_{n-1}}=\frac{2}{\pi}.
\end{equation}

A permutation $w$ \emph{contains} a permutation $\pi\in S_k$ as a \emph{pattern} if $w$ has a subsequence of length $k$ whose entries are in the same relative order as the entries of $\pi$; otherwise $w$ \emph{avoids} $\pi$. We write $\cA_n(\pi)$ and $\cA'_n(\pi)$ for the sets of permutations in $\cA_n$ and $\cA'_n$, respectively, that avoid $\pi$.

Unless otherwise stated, all random permutations and standard Young tableaux are chosen uniformly from the corresponding finite sets. We write $\Pr$ for probability, $\Ex$ for expectation and $\Var$ for variance.

Entringer~\cite{entringer1966combinatorial} refined the Euler numbers by recording the first entry. With the up-down convention used here, let
\begin{equation}\label{eq:entringer-def}
E_{n,k}=\#\{\sigma\in\cA_n:\sigma_1=k\},
\qquad n\ge1,\quad 1\le k\le n.
\end{equation}
These \emph{Entringer numbers} satisfy
\begin{equation}\label{eq:entringer-rec}
E_{n,k}=E_{n,k+1}+E_{n-1,n-k},
\qquad n\ge2,\quad 1\le k<n,
\end{equation}
with initial and boundary conditions $E_{1,1}=1$ and $E_{n,n}=0$ for $n\ge2$. In the OEIS~\cite{OEIS}, the numbers $E_{n,k}$ form the triangle A008281, with each row read from right to left. In particular, $\sum_{k=1}^{n}E_{n,k}=E_n$, and $\Pr(\sigma_1=k)=E_{n,k}/E_n$, where $\sigma$ is chosen uniformly from $\cA_n$. Thus weighted sums of Entringer numbers determine the mean and variance of the first entry.

A \emph{linear extension} of a finite poset $P$ with $N$ elements is an order-preserving bijection from $P$ to $\{1,\ldots,N\}$. We write $L(P)$ for the number of linear extensions of $P$. For $N\ge1$, let $Z_N$ denote the \emph{zigzag poset} on the set $\{1,\ldots,N\}$ of positions, with the relations $1<2>3<4>\cdots$, that is, $2i-1<2i$ and $2i+1<2i$ whenever both sides lie in $\{1,\ldots,N\}$. The linear extensions of $Z_N$ are exactly the up-down alternating permutations $\sigma_1\sigma_2\cdots\sigma_N$ with $\sigma_i=\sigma(i)$, so $L(Z_N)=E_N$. In the figures below, the elements of a poset are drawn as points, with an arrow from $a$ to $b$ when $a<b$ is a cover relation.

A \emph{partition} of an integer $n\geq 0$ is a sequence $\lambda=(\lambda_{1},\lambda_{2},\dots,\lambda_{k})$ of positive integers satisfying $\lambda_{1}\geq\lambda_{2}\geq\cdots\geq\lambda_{k}$ and $\lambda_1+\cdots+\lambda_k=n$. We write $\lambda\vdash n$, $|\lambda|=n$ for the size of $\lambda$, and $\ell(\lambda)=k$ for the number of parts. We also write $n^k=(n,\ldots,n)$ for the partition with $k$ parts equal to $n$. Drawing $\lambda_i$ cells in the $i$-th row, with rows numbered from top to bottom and left justified, gives the \emph{Ferrers diagram} of $\lambda$, which we identify with $\lambda$. If $\lambda$ and $\mu$ are partitions with $\mu\subseteq\lambda$, that is, $\mu_{i}\leq\lambda_{i}$ for all $i$, then the \emph{skew shape} $\lambda/\mu$ is the set of cells of $\lambda$ that are not cells of $\mu$. A \emph{standard Young tableau} of shape $\lambda/\mu$ is a filling of the cells of $\lambda/\mu$ with the integers $1,2,\ldots,|\lambda/\mu|$, each used once, such that the entries increase along each row from left to right and along each column from top to bottom. We write $\SYT(\lambda/\mu)$ for the set of standard Young tableaux of shape $\lambda/\mu$ and $f(\lambda/\mu)$ for their number; when $\mu$ is empty we simply write $\SYT(\lambda)$ and $f(\lambda)$. For $T\in\SYT(\lambda)$ we write $T(i,j)$ for the entry of $T$ in row $i$ and column $j$. The hook length formula of Frame, Robinson and Thrall~\cite{FRT54} states that
\begin{equation}\label{eq:hook}
f(\lambda)=\frac{|\lambda|!}{\prod_{u\in\lambda}h(u)},
\end{equation}
where $h(u)$ is the hook length of the cell $u$, the number of cells of $\lambda$ weakly to the right of $u$ in its row or weakly below $u$ in its column. For skew shapes we use Aitken's determinantal formula
\begin{equation}\label{eq:aitken}
f(\lambda/\mu)
=|\lambda/\mu|!\,
\det\left(
\frac{1}{(\lambda_i-\mu_j-i+j)!}
\right)_{i,j=1}^{\ell(\lambda)},
\end{equation}
where $\mu_j=0$ for $j>\ell(\mu)$ and $1/r!=0$ for negative integers $r$; see Stanley~\cite{ec2} for both formulas. By convention, $\binom{p}{r}=0$ if $r<0$ or $r>p$.

\section{First entry statistics of alternating permutations}

Let $\sigma=\sigma_1\sigma_2\cdots\sigma_{N}$ be chosen uniformly from $\cA_{N}$. In this section, we give two basic statistical results for $\sigma_1$: its expectation $\Ex[\sigma_1]$ and its variance $\Var(\sigma_1)$. Both come from counting linear extensions of the posets that arise from $Z_N$ by adding one or two points.

\subsection{The expectation}
Take $N=4$ as an example. The five permutations in $\cA_4$ are
\[
1324\phantom{1234} 1423\phantom{1234} 2314\phantom{1234} 2413\phantom{1234} 3412.
\]
Thus $\Ex[\sigma_1]=(1+1+2+2+3)/5=1.8$.

\begin{figure}[!ht]
	\centering
	\begin{tikzpicture}[line width=1pt]
		\filldraw (-0.1,-0.1) rectangle (0.1,0.1);
	    \filldraw (0.8,0.8) circle (0.1);
	    \filldraw (1.6,0) circle (0.1);
	    \filldraw (2.4,0.8) circle (0.1);

	    \draw[-latex] (0,0) -- (0.73,0.73);
	    \draw[-latex] (1.6,0) -- (0.87,0.73);
	    \draw[-latex] (1.6,0) -- (2.33,0.73);
	\end{tikzpicture}
	\caption{The zigzag poset $Z_4$ of $\cA_4$; the square is the first position}
	\label{fig:A4}
\end{figure}
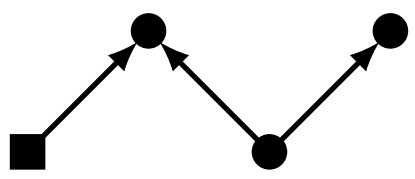

For general $N$, we count linear extensions of posets. We first explain the idea for $N=4$. The number of elements of $\cA_{4}$, namely $E_4$, equals the number of linear extensions of the zigzag poset $Z_4$ shown in Figure~\ref{fig:A4}, where the first position is marked as a square.

We then add one point $v$ above the square, joined to it by an arrow. The new poset $P$ is the zigzag poset of the down-up alternating permutations of length $5$, so $L(P)=E_5=16$.

\begin{figure}[!ht]
	\centering
	\begin{tikzpicture}[line width=1pt]
		\filldraw (0.7,-0.1) rectangle (0.9,0.1);
		\filldraw (0,0.8) circle (0.1);
		\filldraw (1.6,0.8) circle (0.1);
		\filldraw (2.4,0) circle (0.1);
		\filldraw (3.2,0.8) circle (0.1);

		\draw[-latex] (0.8,0) -- (0.07,0.73);
		\draw[-latex] (0.8,0) -- (1.53,0.73);
		\draw[-latex] (2.4,0) -- (1.67,0.73);
		\draw[-latex] (2.4,0) -- (3.13,0.73);

		\draw (3.7,0.35) -- (4,0.35);
		\draw (3.7,0.45) -- (4,0.45);

		\filldraw (5.1,-0.1) rectangle (5.3,0.1);
		\filldraw (4.4,0.8) circle (0.1);
		\filldraw (6,0.8) circle (0.1);
		\filldraw (6.8,0) circle (0.1);
		\filldraw (7.6,0.8) circle (0.1);

		\draw[-latex] (5.2,0) -- (5.93,0.73);
		\draw[-latex] (6.8,0) -- (6.07,0.73);
		\draw[-latex] (6.8,0) -- (7.53,0.73);

		\draw (8,0.4) -- (8.3,0.4);

		\filldraw (9.4,-0.1) rectangle (9.6,0.1);
		\filldraw (8.7,0.8) circle (0.1);
		\filldraw (10.3,0.8) circle (0.1);
		\filldraw (11.1,0) circle (0.1);
		\filldraw (11.9,0.8) circle (0.1);

		\draw[-latex] (8.77,0.73) -- (9.4,0.1);
		\draw[-latex] (9.5,0) -- (10.23,0.73);
		\draw[-latex] (11.1,0) -- (10.37,0.73);
		\draw[-latex] (11.1,0) -- (11.83,0.73);
	\end{tikzpicture}
	\caption{$L(P)=L(P_0)-L(P_1)$ for $N=4$}
	\label{fig:A5eqn}
\end{figure}
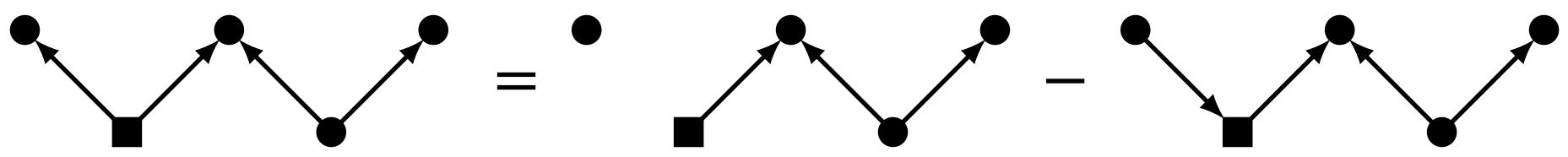

Let $P_0$ be the poset obtained by removing the arrow between $v$ and the square, and let $P_1$ be the poset obtained by reversing this arrow, so that $v$ lies below the square. These posets are shown in Figure~\ref{fig:A5eqn}. A linear extension of $P_0$ puts $v$ either above or below the square, so it is a linear extension of exactly one of $P$ and $P_1$. Hence
\begin{equation}\label{eq:incl-excl}
L(P)=L(P_0)-L(P_1).
\end{equation}

In $P_0$ the point $v$ is incomparable with every element of $Z_4$, so $L(P_0)=\binom{5}{1}E_4$. In $P_1$ the point $v$ is less than the square, so for each $\sigma\in\cA_4$ there are $\sigma_1$ ways to choose the label of $v$: choose $w\in\{1,\ldots,\sigma_1\}$, label $v$ with $w$, and increase every entry of $\sigma$ that is at least $w$ by one. Hence $L(P_1)=\sum_{\sigma\in\cA_4}\sigma_1=\Ex[\sigma_1]E_4$, and \eqref{eq:incl-excl} becomes
\[
E_5=\binom{5}{1}E_4-\Ex[\sigma_1]E_4.
\]
From this we get $\Ex[\sigma_1]=5-E_5/E_4=5-16/5=1.8$, which agrees with the result obtained by direct enumeration.

The same construction works for every $N\ge1$: adding one point above the first position of $Z_N$ gives the zigzag poset of the down-up alternating permutations of length $N+1$, and we get
\begin{equation}\label{eq:mean-identity}
E_{N+1}=(N+1)E_N-\sum_{\sigma\in\cA_N}\sigma_1.
\end{equation}

\subsection{The variance}
Using a similar method, we can determine the variance of $\sigma_1$. We extend $Z_4$ to $Z_6$ by adding two points $v_1<v_2$ to the left of the square, with $v_2$ above the square.

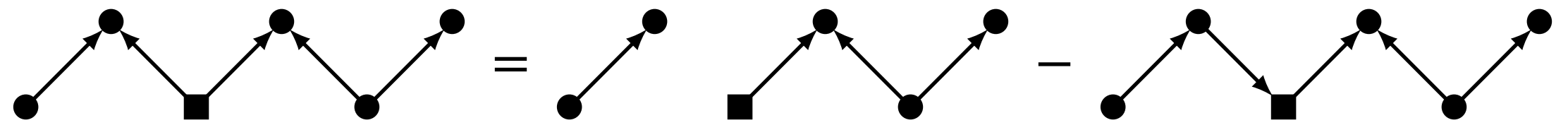
\begin{figure}[!ht]
	\centering
	\begin{tikzpicture}[line width=1pt]
		\filldraw (0,0) circle (0.1);
		\filldraw (1.5,-0.1) rectangle (1.7,0.1);
		\filldraw (0.8,0.8) circle (0.1);
		\filldraw (2.4,0.8) circle (0.1);
		\filldraw (3.2,0) circle (0.1);
		\filldraw (4,0.8) circle (0.1);

		\draw[-latex] (0,0) -- (0.73,0.73);
		\draw[-latex] (1.6,0) -- (0.87,0.73);
		\draw[-latex] (1.6,0) -- (2.33,0.73);
		\draw[-latex] (3.2,0) -- (2.47,0.73);
		\draw[-latex] (3.2,0) -- (3.93,0.73);

		\draw (4.4,0.35) -- (4.7,0.35);
		\draw (4.4,0.45) -- (4.7,0.45);

		\filldraw (5.1,0) circle (0.1);
		\filldraw (6.6,-0.1) rectangle (6.8,0.1);
		\filldraw (5.9,0.8) circle (0.1);
		\filldraw (7.5,0.8) circle (0.1);
		\filldraw (8.3,0) circle (0.1);
		\filldraw (9.1,0.8) circle (0.1);

		\draw[-latex] (5.1,0) -- (5.83,0.73);
		\draw[-latex] (6.7,0) -- (7.43,0.73);
		\draw[-latex] (8.3,0) -- (7.57,0.73);
		\draw[-latex] (8.3,0) -- (9.03,0.73);

		\draw (9.5,0.4) -- (9.8,0.4);

		\filldraw (10.2,0) circle (0.1);
		\filldraw (11.7,-0.1) rectangle (11.9,0.1);
		\filldraw (11,0.8) circle (0.1);
		\filldraw (12.6,0.8) circle (0.1);
		\filldraw (13.4,0) circle (0.1);
		\filldraw (14.2,0.8) circle (0.1);

		\draw[-latex] (10.2,0) -- (10.93,0.73);
		\draw[-latex] (11.07,0.73) -- (11.7,0.1);
		\draw[-latex] (11.8,0) -- (12.53,0.73);
		\draw[-latex] (13.4,0) -- (12.67,0.73);
		\draw[-latex] (13.4,0) -- (14.13,0.73);
	\end{tikzpicture}
	\caption{$L(P)=L(P_0)-L(P_1)$ for $N=4$}
	\label{fig:A6eqn}
\end{figure}

Let $P=Z_6$, let $P_0$ be the poset obtained by removing the arrow between $v_2$ and the square, and let $P_1$ be the poset obtained by reversing this arrow. These posets are shown in Figure~\ref{fig:A6eqn}, and again $L(P)=L(P_0)-L(P_1)$.

In $P_0$, each element of the chain $v_1<v_2$ is incomparable with every element of $Z_4$. Thus, choosing the two labels assigned to the chain and a linear extension of $Z_4$ gives $L(P_0)=\binom{6}{2}E_4$. In $P_1$ the two new points are less than the square, so for each $\sigma\in\cA_4$ there are $\binom{\sigma_1+1}{2}$ ways to choose their labels, namely a $2$-subset of $\{1,\ldots,\sigma_1+1\}$. Then we can rewrite $L(P)=L(P_0)-L(P_1)$ as
\[
E_6=\binom{6}{2} E_4-\Ex\left[\binom{\sigma_1+1}{2}\right]E_4,
\]
and hence we get $\Ex[\sigma_1^2+\sigma_1]=30-2E_6/E_4=30-122/5=5.6$. For general $N$, adding two points in the same way turns $Z_N$ into $Z_{N+2}$, and we get
\begin{equation}\label{eq:second-identity}
E_{N+2}=\binom{N+2}{2}E_N-\sum_{\sigma\in\cA_N}\binom{\sigma_1+1}{2}.
\end{equation}

\begin{theorem}\label{thm:first-entry}
Let $N\ge1$ and let $\sigma=\sigma_1\sigma_2\cdots\sigma_N$ be chosen uniformly from $\cA_N$. Then
\begin{equation}\label{eq:mean}
\Ex[\sigma_1]=N+1-\frac{E_{N+1}}{E_{N}},
\end{equation}
\begin{equation}\label{eq:binom-moment}
\Ex\left[\binom{\sigma_1+1}{2}\right]=\binom{N+2}{2}-\frac{E_{N+2}}{E_{N}},
\end{equation}
and
\begin{equation}\label{eq:var}
\Var(\sigma_1)=\frac{(2N+3)E_{N}E_{N+1}-E^2_{N+1}-2E_{N}E_{N+2}}{E^2_{N}}.
\end{equation}
Equivalently, in terms of Entringer numbers,
\begin{equation}\label{eq:entringer-sums}
\sum_{k=1}^{N}k\,E_{N,k}=(N+1)E_N-E_{N+1},
\qquad
\sum_{k=1}^{N}\binom{k+1}{2}E_{N,k}=\binom{N+2}{2}E_N-E_{N+2}.
\end{equation}
\end{theorem}

\begin{proof}
Dividing \eqref{eq:mean-identity} and \eqref{eq:second-identity} by $E_N$ gives \eqref{eq:mean} and \eqref{eq:binom-moment}. Since $\sigma_1^2=2\binom{\sigma_1+1}{2}-\sigma_1$, we have
\[
\Ex[\sigma_1^2]=2\Ex\left[\binom{\sigma_1+1}{2}\right]-\Ex[\sigma_1]=(N+1)^2-\frac{2E_{N+2}-E_{N+1}}{E_{N}}.
\]
Therefore
\begin{align*}
\Var(\sigma_1)
&=\Ex[\sigma_1^2]-\Ex[\sigma_1]^2\\
&=(N+1)^2-\frac{2E_{N+2}-E_{N+1}}{E_{N}}-\left(N+1-\frac{E_{N+1}}{E_{N}}\right)^2\\
&=\frac{(2N+3)E_{N+1}}{E_N}-\frac{2E_{N+2}}{E_N}-\frac{E_{N+1}^2}{E_N^2},
\end{align*}
which is \eqref{eq:var}. Finally, since $\sum_{\sigma\in\cA_N}g(\sigma_1)=\sum_{k=1}^{N}g(k)E_{N,k}$ for every function $g$, the identities \eqref{eq:entringer-sums} are just \eqref{eq:mean-identity} and \eqref{eq:second-identity} written in terms of Entringer numbers.
\end{proof}

For $N=4$, formula \eqref{eq:var} gives $\Var(\sigma_1)=14/25=0.56$, which agrees with the direct computation $\Ex[\sigma_1^2]-\Ex[\sigma_1]^2=3.8-3.24=0.56$. The identities \eqref{eq:entringer-sums} can also be derived directly from the recurrence \eqref{eq:entringer-rec}, which gives
\[
E_{n,k}=\sum_{m=1}^{n-k}E_{n-1,m}.
\]

When $N$ is large, the asymptotic behavior \eqref{eq:ratio} of the Euler numbers gives the following result.

\begin{corollary}\label{cor:first-entry-asym}
As $N\to\infty$,
\[
\Ex[\sigma_1]=\left(1-\frac{2}{\pi}\right)(N+1)+O\left(N3^{-N}\right),
\qquad
\Var(\sigma_1)=\frac{4(\pi-3)}{\pi^2}\,N^2+O(N).
\]
In particular, $\Ex[\sigma_1]\sim(1-2/\pi)N\approx 0.36N$ and $\Var(\sigma_1)\sim4(\pi-3)N^2/\pi^2\approx0.057N^2$.
\end{corollary}

\begin{proof}
By \eqref{eq:ratio},
\[
\frac{E_{N+1}}{E_N}=\frac{2}{\pi}(N+1)\bigl(1+O(3^{-N})\bigr),
\qquad
\frac{E_{N+2}}{E_N}=\frac{4}{\pi^2}(N+1)(N+2)\bigl(1+O(3^{-N})\bigr).
\]
Substituting the first estimate into \eqref{eq:mean} gives the formula for $\Ex[\sigma_1]$. Substituting both estimates into \eqref{eq:var} gives
\begin{align*}
\Var(\sigma_1)
&=\frac{2}{\pi}(2N+3)(N+1)-\frac{8}{\pi^2}(N+1)(N+2)-\frac{4}{\pi^2}(N+1)^2+O\left(N^23^{-N}\right)\\
&=\left(\frac{4}{\pi}-\frac{12}{\pi^2}\right)N^2+O(N)\\
&=\frac{4(\pi-3)}{\pi^2}\,N^2+O(N).
\qedhere
\end{align*}
\end{proof}

The constants $1-2/\pi$ and $4(\pi-3)/\pi^2$ in Corollary~\ref{cor:first-entry-asym} are the mean and the variance of a limit law of Simons and Yao~\cite{simons1996sine}. Their second sine law, transferred to the up-down convention by the complement map, states that for $\sigma$ chosen uniformly from $\cA_N$,
\[
\Pr(\sigma_1=k)=\frac{\pi}{2N}\cos\frac{(2k-1)\pi}{4N}+O\left(N^{-2}\right)
\]
uniformly in $1\le k\le N$. Thus $\sigma_1/N$ converges in distribution to the law on $[0,1]$ with density $(\pi/2)\cos(\pi x/2)$. Since
\[
\int_0^1 x\,\frac{\pi}{2}\cos\frac{\pi x}{2}\,dx=1-\frac{2}{\pi},
\qquad
\int_0^1 x^2\,\frac{\pi}{2}\cos\frac{\pi x}{2}\,dx=1-\frac{8}{\pi^2},
\]
this law has mean $1-2/\pi$ and variance $1-8/\pi^2-(1-2/\pi)^2=4(\pi-3)/\pi^2$, in agreement with Corollary~\ref{cor:first-entry-asym}.

By the complement bijection between $\cA_N$ and $\cA'_N$, the results on the first entry carry over to down-up alternating permutations.

\begin{corollary}\label{cor:down-up-first}
Let $\sigma'=\sigma'_1\sigma'_2\cdots\sigma'_N$ be chosen uniformly from $\cA'_N$. Then
\[
\Ex[\sigma'_1]=N+1-\Ex[\sigma_1]=\frac{E_{N+1}}{E_N},
\]
and $\Var(\sigma'_1)=\Var(\sigma_1)$ is given by \eqref{eq:var}. In particular, $\Ex[\sigma'_1]=2(N+1)/\pi+O(N3^{-N})$.
\end{corollary}

\begin{proof}
Under the complement map, $\sigma'_1=N+1-\sigma_1$. Hence $\sigma'_1$ has the same distribution as $N+1-\sigma_1$, and the claims follow from Theorem~\ref{thm:first-entry} and Corollary~\ref{cor:first-entry-asym}.
\end{proof}

\section{Positions of fixed values in alternating permutations}

In this section, we study the positions of fixed values in random alternating permutations. We first use a bivariate exponential generating function to compute the expectation and the variance of the position of the value $1$. We then obtain the corresponding formulas for every other value by interchanging consecutive values. The formulas depend on the parity of the length.

For $\sigma\in\cA_N$ and $1\le j\le N$, let
\begin{equation}\label{eq:Xj}
X_j(\sigma)=\sigma^{-1}(j)
\end{equation}
denote the position of the value $j$ in $\sigma$.

\begin{figure}[!ht]
	\centering
	\begin{tikzpicture}[line width=1pt]
		\filldraw (0,0) circle (0.1);
		\filldraw (0.8,0.8) circle (0.1);
		\filldraw (1.6,0) circle (0.1);
		\filldraw (2.4,0.8) circle (0.1);
		\filldraw (4.4,0.8) circle (0.1);
		\filldraw (5.113,0.05) -- (5.286,0.05) -- (5.2,-0.1) -- cycle;
		\filldraw (6,0.8) circle (0.1);
		\filldraw (6.8,0) circle (0.1);
		\filldraw (7.6,0.8) circle (0.1);
		\filldraw (9.6,0) circle (0.1);
		\filldraw (10.4,0.8) circle (0.1);

		\draw[-latex] (0,0) -- (0.73,0.73);
		\draw[-latex] (1.6,0) -- (0.87,0.73);
		\draw[-latex] (1.6,0) -- (2.33,0.73);
		\draw[-latex] (5.2,0) -- (4.47,0.73);
		\draw[-latex] (5.2,0) -- (5.93,0.73);
		\draw[-latex] (6.8,0) -- (6.07,0.73);
		\draw[-latex] (6.8,0) -- (7.53,0.73);
		\draw[-latex] (9.6,0) -- (10.33,0.73);

		\node at (3.4,0.4) {$\cdots$};
		\node at (8.6,0.4) {$\cdots$};

		\draw[decorate,decoration={brace,amplitude=2mm,mirror,raise=5pt}] (0,-0.15) -- (4.4,-0.15);
		\draw[->] (5.2,-0.2) -- (5.2,-0.6);
		\draw[decorate,decoration={brace,amplitude=2mm,mirror,raise=5pt}] (6,-0.15) -- (10.4,-0.15);

		\node at (2.2,-0.85) {$2m$};
		\node at (5.2,-0.85) {$1$};
		\node at (8.2,-0.85) {$2n-1-2m$};
	\end{tikzpicture}
	\caption{The decomposition of $\sigma\in\cA_{2n}$ when the entry $1$ occurs at position $2m+1$}
	\label{fig:one-position-even}
\end{figure}
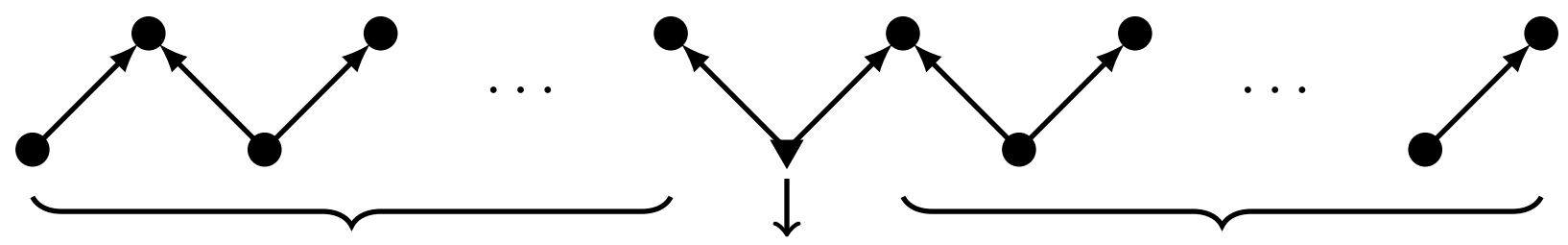

We first consider the position of the minimum value. Since $1$ must be a valley in an up-down alternating permutation, it can occur only in an odd position. If $X_1=2m+1$, then the entries to the left of $1$ form an alternating block of even length $2m$, while the entries to the right form an alternating block of length $N-1-2m$, as shown in Figure~\ref{fig:one-position-even} for even $N$. After choosing the values used on the left, the number of such permutations is
\begin{equation}\label{eq:a}
a_{N,2m+1}
=\binom{N-1}{2m}E_{2m}E_{N-1-2m},
\qquad 0\le m\le\left\lfloor\frac{N-1}{2}\right\rfloor.
\end{equation}
Equivalently,
\[
\Pr(X_1=2m+1)
=\frac{\binom{N-1}{2m}E_{2m}E_{N-1-2m}}{E_{N}}.
\]

To compute the moments of $X_1$, we use the bivariate exponential generating function
\[
F(x,u)
=
\sum_{N\ge1}
\left(
\sum_{\sigma\in\cA_{N}}
u^{X_1(\sigma)}
\right)
\frac{x^{N-1}}{(N-1)!}.
\]
The variable $x$ marks the entries greater than $1$, and the variable $u$ marks the position of $1$. Every entry of the left block contributes one unit to $X_1$, so the left block is generated by $\sec(ux)$. The entries of the right block do not affect $X_1$, so by \eqref{eq:andre} the right block is generated by $\sec x+\tan x$. The value $1$ itself contributes the factor $u$. Hence
\begin{equation}\label{eq:F1}
F(x,u)=u\sec(ux)\bigl(\sec x+\tan x\bigr).
\end{equation}
Expanding,
\begin{align*}
F(x,u)
&=u\left(\sum_{m\geq 0}E_{2m}u^{2m}\frac{x^{2m}}{(2m)!}\right)
\left(\sum_{r\geq 0}E_{r}\frac{x^{r}}{r!}\right)\\
&=\sum_{m,r\geq 0}
E_{2m}E_{r}u^{2m+1}
\frac{x^{2m+r}}{(2m)!\,r!},
\end{align*}
so the coefficient of $x^{N-1}/(N-1)!$ is $\sum_{m}a_{N,2m+1}u^{2m+1}$, and setting $u=1$ gives the total number $E_{N}$. Since $\sec(ux)$ and $\sec x$ are even functions of $x$ and $\tan x$ is odd, the terms with $N$ even form $u\sec(ux)\tan x$, and the terms with $N$ odd form $u\sec(ux)\sec x$.

\begin{proposition}\label{prop:min-position}
Let $N\ge2$, let $\sigma$ be chosen uniformly from $\cA_N$, and put $\rho=E_{N-2}/E_N$.
\begin{enumerate}[(a)]
\item If $N=2n$, then
\begin{equation}\label{eq:EX1}
\Ex[X_1]=n+\frac{1}{2}-\frac{2n-1}{2}\,\rho
\end{equation}
and
\begin{equation}\label{eq:VarX1even}
\Var(X_1)=\frac{4n^2-1}{12}-\frac{(n-1)(2n-1)}{3}\,\rho-\frac{(2n-1)^2}{4}\,\rho^2.
\end{equation}
\item If $N=2n+1$, then $\Ex[X_1]=n+1$ and
\begin{equation}\label{eq:VarX1odd}
\Var(X_1)=\frac{n(n+1)}{3}+\frac{2n(2n-1)}{3}\,\rho.
\end{equation}
\end{enumerate}
\end{proposition}

\begin{proof}
In both cases $\Ex[X_1]=\sum_{m}(2m+1)\,a_{N,2m+1}/E_{N}$, so we differentiate the corresponding part of \eqref{eq:F1} with respect to $u$ and then set $u=1$; in the same way, $\sum_{\sigma\in\cA_N}X_1(\sigma)\bigl(X_1(\sigma)-1\bigr)$ is the coefficient of $x^{N-1}/(N-1)!$ in the second derivative at $u=1$. We then use
\begin{equation}\label{eq:var-from-fact}
\Var(X_1)=\Ex[X_1(X_1-1)]+\Ex[X_1]-\Ex[X_1]^2.
\end{equation}

(a) Let $N=2n$. Write $s(x)=\sec x$. We will use $s'(x)=\sec x\tan x$, $s''(x)=2\sec^3x-\sec x$ and $s'''(x)=6\sec^3x\tan x-\sec x\tan x$, so that
\begin{equation}\label{eq:sec-identities}
s''(x)-s(x)=2\sec x\tan^2x,
\qquad
s'''(x)-2s'(x)=3\sec x\tan x\,(2\sec^2x-1).
\end{equation}
The first derivative is
\begin{align*}
\left.\frac{\partial}{\partial u}\bigl(u\sec(ux)\tan x\bigr)\right|_{u=1}
&=\sec x\tan x+x\sec x\tan^2x\\
&=s'(x)+\frac{x}{2}\bigl(s''(x)-s(x)\bigr).
\end{align*}
By \eqref{eq:sectan}, the coefficients of $x^{2n-1}/(2n-1)!$ in $s'(x)$, $xs''(x)$ and $xs(x)$ are
\[
E_{2n},\qquad
(2n-1)E_{2n},\qquad
(2n-1)E_{2n-2},
\]
respectively. Therefore
\begin{align*}
\Ex[X_1]
&=\frac{1}{E_{2n}}
\left(E_{2n}+\frac{2n-1}{2}E_{2n}-\frac{2n-1}{2}E_{2n-2}\right)\\
&=n+\frac{1}{2}-\frac{2n-1}{2}\,\rho.
\end{align*}
The second derivative is, by \eqref{eq:sec-identities},
\begin{align*}
\left.\frac{\partial^2}{\partial u^2}\bigl(u\sec(ux)\tan x\bigr)\right|_{u=1}
&=2x\sec x\tan^2x+x^2\sec x\tan x\,(2\sec^2x-1)\\
&=x\bigl(s''(x)-s(x)\bigr)
+
\frac{x^2}{3}\bigl(s'''(x)-2s'(x)\bigr).
\end{align*}
By \eqref{eq:sectan}, the coefficients of $x^{2n-1}/(2n-1)!$ in $x^2s'''(x)$ and $x^2s'(x)$ are $(2n-1)(2n-2)E_{2n}$ and $(2n-1)(2n-2)E_{2n-2}$. Hence
\begin{align*}
\sum_{\sigma\in\cA_{2n}}
X_1(\sigma)\bigl(X_1(\sigma)-1\bigr)
&=(2n-1)\bigl(E_{2n}-E_{2n-2}\bigr)+\frac{(2n-1)(2n-2)}{3}\bigl(E_{2n}-2E_{2n-2}\bigr)\\
&=\frac{4n^2-1}{3}E_{2n}
+
\frac{-8n^2+6n-1}{3}E_{2n-2}.
\end{align*}
We divide by $E_{2n}$ and substitute into \eqref{eq:var-from-fact}. In the result, the constant term and the coefficient of $\rho$ are
\begin{align*}
\frac{4n^2-1}{3}+n+\frac12-\left(n+\frac12\right)^2&=\frac{4n^2-1}{12},\\
\frac{-8n^2+6n-1}{3}-\frac{2n-1}{2}+\frac{(2n+1)(2n-1)}{2}&=-\frac{(n-1)(2n-1)}{3},
\end{align*}
and the coefficient of $\rho^2$ is $-(2n-1)^2/4$. This proves \eqref{eq:VarX1even}.

(b) Let $N=2n+1$. Reversal maps $\delta=\delta_1\cdots\delta_{2n+1}\in\cA_{2n+1}$ to $\delta_{2n+1}\cdots\delta_1$, which is again up-down alternating because $2n+1$ is odd. If a value is at position $k$ in $\delta$, then it is at position $2n+2-k$ after reversal. Hence the distribution of $X_1$ is symmetric about $n+1$, and $\Ex[X_1]=n+1$. Using $(\tan x)''=2\sec^2x\tan x$ and $(\tan x)'''=6\sec^4x-4\sec^2x$, the second derivative is
\begin{align*}
\left.\frac{\partial^2}{\partial u^2}\bigl(u\sec(ux)\sec x\bigr)\right|_{u=1}
&=2x\sec^2x\tan x+x^2\sec^2x\,(2\sec^2x-1)\\
&=x\,(\tan x)''+\frac{x^2}{3}\bigl((\tan x)'''+(\tan x)'\bigr).
\end{align*}
By \eqref{eq:sectan}, the coefficients of $x^{2n}/(2n)!$ in $x(\tan x)''$, $x^2(\tan x)'''$ and $x^2(\tan x)'$ are
\[
2nE_{2n+1},\qquad
2n(2n-1)E_{2n+1},\qquad
2n(2n-1)E_{2n-1},
\]
respectively. Hence
\[
\sum_{\sigma\in\cA_{2n+1}}
X_1(\sigma)\bigl(X_1(\sigma)-1\bigr)
=2nE_{2n+1}+\frac{2n(2n-1)}{3}\bigl(E_{2n+1}+E_{2n-1}\bigr),
\]
and by \eqref{eq:var-from-fact},
\[
\Var(X_1)
=2n+\frac{2n(2n-1)}{3}(1+\rho)+(n+1)-(n+1)^2
=\frac{n(n+1)}{3}+\frac{2n(2n-1)}{3}\,\rho.
\qedhere
\]
\end{proof}

For $j\ge2$, one can decompose $\sigma$ in the same way around the values $1,\ldots,j$, but the number of cases grows with $j$. Instead, we pass from the value $1$ to an arbitrary value $j$ by comparing the positions of $j$ and $j+1$.

\begin{lemma}\label{lem:position-sum-recurrences}
For $N\ge2$ and $1\le j<N$,
\begin{equation}\label{eq:shift1}
\sum_{\sigma\in\cA_N}
\bigl(\sigma^{-1}(j+1)-\sigma^{-1}(j)\bigr)
=
\begin{cases}
E_{N-2}, & \text{if $N$ is even},\\
0, & \text{if $N$ is odd},
\end{cases}
\end{equation}
and
\begin{equation}\label{eq:shift2}
\sum_{\sigma\in\cA_N}
\left(
\bigl(\sigma^{-1}(j+1)\bigr)^2
-
\bigl(\sigma^{-1}(j)\bigr)^2
\right)
=
\begin{cases}
(2N+1-2j)E_{N-2}, & \text{if $N$ is even},\\
-2(j-1)E_{N-2}, & \text{if $N$ is odd}.
\end{cases}
\end{equation}
\end{lemma}

\begin{proof}
For $N=2$, the only permutation is $12$, and both identities follow directly. We therefore assume $N\ge3$ throughout the remainder of the proof.
If $j$ and $j+1$ are not adjacent in $\sigma$, interchanging them preserves all alternating inequalities. This involution changes the signs of the differences of their positions and of their squared positions. Thus these contributions cancel in both sums.

If $j$ and $j+1$ are adjacent, deleting the pair and standardizing gives an alternating permutation of length $N-2$. To see this, note that an increasing pair $j,j+1$ is preceded, when present, by a value greater than $j+1$ and followed, when present, by a value smaller than $j$. For a decreasing pair $j+1,j$, these boundary inequalities are reversed. Hence deletion preserves the alternating pattern.

Fix $\sigma\in\cA_{N-2}$. To recover the original labels, replace every entry at least $j$ by that entry plus $2$. Consider inserting the pair into the gap after the first $t$ entries, where $0\le t\le N-2$. For convenience, set $\sigma_0=N$, and set $\sigma_{N-1}=0$ if $N$ is even and $\sigma_{N-1}=N$ if $N$ is odd. For a statement $A$, let $\chi(A)=1$ if $A$ is true and $\chi(A)=0$ if $A$ is false.

For even $t$, an increasing pair $j,j+1$ can be inserted exactly when $\sigma_t\ge j$ and $\sigma_{t+1}<j$. For odd $t$, a decreasing pair $j+1,j$ can be inserted exactly when $\sigma_t<j$ and $\sigma_{t+1}\ge j$. By the choice of $\sigma_0$ and $\sigma_{N-1}$, we have $\sigma_t>\sigma_{t+1}$ for even $t$ and $\sigma_t<\sigma_{t+1}$ for odd $t$, also at the two ends. Hence the signed contribution of this gap to the sum in \eqref{eq:shift1} is
\[
\chi(\sigma_{t+1}<j)
-
\chi(\sigma_t<j).
\]
Summing over all gaps gives
\[
\sum_{t=0}^{N-2}
\left(
\chi(\sigma_{t+1}<j)
-
\chi(\sigma_t<j)
\right)
=\chi(\sigma_{N-1}<j)-\chi(\sigma_{0}<j)
=\chi(\sigma_{N-1}<j),
\]
which is $1$ if $N$ is even and $0$ if $N$ is odd. Summing over $\sigma\in\cA_{N-2}$ proves \eqref{eq:shift1}.

The inserted pair occupies positions $t+1$ and $t+2$. Its signed contribution to the sum in \eqref{eq:shift2} is therefore
\[
(2t+3)
\left(
\chi(\sigma_{t+1}<j)
-
\chi(\sigma_t<j)
\right).
\]
Summing over the gaps and using that exactly $j-1$ entries of $\sigma$ are smaller than $j$ gives
\begin{align*}
&\sum_{t=0}^{N-2}(2t+3)
\left(
\chi(\sigma_{t+1}<j)
-
\chi(\sigma_t<j)
\right)\\
&\qquad=
(2N-1)\chi(\sigma_{N-1}<j)-3\chi(\sigma_{0}<j)-2\sum_{t=1}^{N-2}\chi(\sigma_t<j)\\
&\qquad=
(2N-1)\chi(\sigma_{N-1}<j)-2(j-1),
\end{align*}
which is $2N+1-2j$ if $N$ is even and $-2(j-1)$ if $N$ is odd. Summing over $\sigma\in\cA_{N-2}$ proves \eqref{eq:shift2}.
\end{proof}

\begin{theorem}\label{thm:general-position}
Let $N\ge2$, let $\sigma$ be chosen uniformly from $\cA_N$, let $1\le j\le N$, and put $\rho=E_{N-2}/E_N$.
\begin{enumerate}[(a)]
\item If $N=2n$, then
\begin{equation}\label{eq:EXj}
\Ex[X_j]
=
n+\frac12-\frac{2n-2j+1}{2}\,\rho
=
\frac{2n+1}{2}
+
\left(j-\frac{2n+1}{2}\right)\rho
\end{equation}
and
\begin{equation}\label{eq:VarXj}
\begin{split}
\Var(X_j)
&=
\frac{4n^2-1}{12}
+
\frac{-2n^2+(6j-3)n-3j(j-1)-1}{3}\,\rho\\
&\qquad
-
\frac{(2n-2j+1)^2}{4}\,\rho^2.
\end{split}
\end{equation}
\item If $N=2n+1$, then $\Ex[X_j]=n+1$ and
\begin{equation}\label{eq:VarXjodd}
\Var(X_j)=\frac{n(n+1)}{3}+\left(\frac{2n(2n-1)}{3}-(j-1)(j-2)\right)\rho.
\end{equation}
\end{enumerate}
In particular, for every fixed $j$, $\Ex[X_j]=(N+1)/2+O(N^{-1})$ and $\Var(X_j)=(N^2-1)/12+O(1)$.
\end{theorem}

\begin{proof}
Applying Lemma~\ref{lem:position-sum-recurrences} to the values $1,2,\ldots,j-1$ and summing gives
\begin{equation}\label{eq:shift-sum}
\Ex[X_j]-\Ex[X_1]=
\begin{cases}
(j-1)\rho, & \text{if $N$ is even},\\
0, & \text{if $N$ is odd},
\end{cases}
\end{equation}
and
\[
\Ex[X_j^2]-\Ex[X_1^2]=
\begin{cases}
\displaystyle\sum_{r=1}^{j-1}(2N+1-2r)\,\rho=(j-1)(2N+1-j)\,\rho, & \text{if $N$ is even},\\[3ex]
\displaystyle-\sum_{r=1}^{j-1}2(r-1)\,\rho=-(j-1)(j-2)\,\rho, & \text{if $N$ is odd}.
\end{cases}
\]
We also use
\[
\Var(X_j)-\Var(X_1)=\bigl(\Ex[X_j^2]-\Ex[X_1^2]\bigr)-\bigl(\Ex[X_j]-\Ex[X_1]\bigr)\bigl(\Ex[X_j]+\Ex[X_1]\bigr).
\]

(a) Let $N=2n$. By \eqref{eq:EX1} and \eqref{eq:shift-sum},
\[
\Ex[X_j]
=
n+\frac12-\frac{2n-1}{2}\rho+(j-1)\rho
=
n+\frac12-\frac{2n-2j+1}{2}\rho,
\]
which is \eqref{eq:EXj}. Hence $\Ex[X_j]+\Ex[X_1]=2n+1-(2n-j)\rho$, and
\begin{align*}
\Var(X_j)-\Var(X_1)
&=(j-1)(4n+1-j)\rho-(j-1)\rho\bigl(2n+1-(2n-j)\rho\bigr)\\
&=(j-1)(2n-j)\,\rho\,(1+\rho).
\end{align*}
Adding \eqref{eq:VarX1even} and using
\begin{align*}
-\frac{(n-1)(2n-1)}{3}+(j-1)(2n-j)&=\frac{-2n^2+(6j-3)n-3j(j-1)-1}{3},\\
-\frac{(2n-1)^2}{4}+(j-1)(2n-j)&=-\frac{(2n-2j+1)^2}{4},
\end{align*}
we obtain \eqref{eq:VarXj}.

(b) Let $N=2n+1$. By \eqref{eq:shift-sum}, $\Ex[X_j]=\Ex[X_1]=n+1$; this also follows from the reversal argument in the proof of Proposition~\ref{prop:min-position}(b), which applies to every value. Then $\Var(X_j)=\Var(X_1)-(j-1)(j-2)\rho$, and \eqref{eq:VarX1odd} gives \eqref{eq:VarXjodd}.

By \eqref{eq:ratio}, $\rho\sim\pi^2/(4N^2)$. Hence the corrections in \eqref{eq:EXj} are $O(N^{-1})$ and the corrections in \eqref{eq:VarXj} and \eqref{eq:VarXjodd} are $O(1)$ for fixed $j$, while the leading terms $(4n^2-1)/12$ for $N=2n$ and $n(n+1)/3$ for $N=2n+1$ are both equal to $(N^2-1)/12$, the variance of a uniformly random position in $\{1,\ldots,N\}$.
\end{proof}

Finally, the complement map $\sigma_i\mapsto N+1-\sigma_i$ changes up-down alternating permutations into down-up ones and exchanges the values $j$ and $N+1-j$. This gives the corresponding results for down-up alternating permutations, and in particular for the position of the maximum value.

\begin{corollary}\label{cor:complement-positions}
Let $\sigma'$ be chosen uniformly from $\cA'_N$ and let $1\le j\le N$. Then the position of $j$ in $\sigma'$ has the same distribution as $X_{N+1-j}$ on $\cA_N$. In particular, for $N=2n$,
\[
\Ex[(\sigma')^{-1}(1)]=\Ex[X_{2n}]
=n+\frac{1}{2}+\frac{2n-1}{2}\,\frac{E_{2n-2}}{E_{2n}}
\]
and
\[
\Ex[(\sigma')^{-1}(2n)]=\Ex[X_{1}]
=n+\frac{1}{2}-\frac{2n-1}{2}\,\frac{E_{2n-2}}{E_{2n}},
\]
both of which are asymptotic to $n+1/2$. For $N=2n+1$, every value has expected position $n+1$ in $\cA'_{2n+1}$ as well.
\end{corollary}

\begin{proof}
The complement map is a bijection from $\cA_N$ to $\cA'_N$ which exchanges the values $j$ and $N+1-j$ and preserves positions. The formulas follow from Theorem~\ref{thm:general-position}.
\end{proof}

\section{Entry statistics in standard Young tableaux}

We now apply the linear extension approach to entry statistics in standard Young tableaux. Recall that a standard Young tableau of shape $\lambda$ can be regarded as a linear extension of the poset of cells of $\lambda$, in which each cell is less than its right neighbor and its lower neighbor. In a standard Young tableau, the upper-left entry is always $1$. We study the entry in the lower-left corner; transposition identifies this statistic with the upper-right entry in the conjugate shape. We begin with the shape $(n,n)$. An insertion construction expresses the sum of the lower-left entries as the number of standard Young tableaux of a skew shape. We then extend this construction to arbitrary shapes and use the hook length formula to obtain explicit expressions in special cases. Next, we determine the distribution of the entry in each cell of the second row of a two-row tableau, and for the shape $(n,n)$ the sum of the entries in each cell of the second row. At the end of the section, we apply the results for the shape $(n,n)$ to $321$-avoiding down-up permutations.

\subsection{The lower-left entry}\label{subsec:lower-left-entry}

Let $\lambda$ be a partition with $k\ge2$ parts, and let $T$ be chosen uniformly from $\SYT(\lambda)$. Riegler and Neumann~\cite{riegler2016playing} obtained the following formula for the mean of the entry in the cell $(2,1)$:
\begin{equation}\label{eq:RN}
\Ex[T(2,1)]=\prod_{i=1}^{k}\frac{|\lambda|+i}{|\lambda|+i-\lambda_i}.
\end{equation}
Proposition~\ref{prop:two-row-corner-sum} and Corollary~\ref{cor:two-row-mean} below are the cases $\lambda=(n,n)$ and $\lambda=(m,n)$ of \eqref{eq:RN}. The insertion construction gives a short direct proof of these cases, and it applies equally to the lower-left entry $T(k,1)$ of a shape with $k$ rows (Theorem~\ref{thm:general-shape}), which for $k\ge3$ is not covered by \eqref{eq:RN}.

For a two-row shape $(m,n)$ with $m\ge n\ge1$, the hook lengths are $m+1,m,\ldots,m-n+2$ and $m-n,\ldots,1$ in the first row (the second list is empty when $m=n$) and $n,\ldots,1$ in the second row, so \eqref{eq:hook} gives
\begin{equation}\label{eq:two-row-f}
f(m,n)
=\frac{(m+n)!\,(m-n+1)}{(m+1)!\,n!}
=\frac{m-n+1}{m+1}\binom{m+n}{n}
=\binom{m+n}{n}-\binom{m+n}{n-1}.
\end{equation}
In particular, the number of standard Young tableaux of shape $(n,n)$ is the $n$-th Catalan number $C_n=\binom{2n}{n}/(n+1)$; see sequence A000108 in the OEIS~\cite{OEIS}. For example, there are $2$ standard Young tableaux of shape $(2,2)$ and $5$ standard Young tableaux of shape $(3,3)$, as shown in Figure~\ref{fig:T3}, where the lower-left entries are marked in boldface.

\begin{figure}[!ht]
	\centering
	\begin{minipage}[b]{0.2\textwidth}
	\centering
		\begin{tikzpicture}[line width=1pt]
			\draw[step=0.5] (0,0 ) grid ( 1, 1);

			\node at (0.25,0.75) { 1 };
			\node at (0.25,0.25) {\textbf{2}};
			\node at (0.75,0.75) { 3 };
			\node at (0.75,0.25) { 4 };
		\end{tikzpicture}
	\end{minipage}
	\begin{minipage}[b]{0.2\textwidth}
		\centering
		\begin{tikzpicture}[line width=1pt]
			\draw[step=0.5] (0,0 ) grid ( 1, 1);

			\node at (0.25,0.75) { 1 };
			\node at (0.75,0.75) { 2 };
			\node at (0.25,0.25) {\textbf{3}};
			\node at (0.75,0.25) { 4 };
		\end{tikzpicture}
	\end{minipage}

		\vspace{2ex}

	\begin{minipage}[b]{0.15\textwidth}
		\centering
		\begin{tikzpicture}[line width=1pt]
			\draw[step=0.5] (0,0 ) grid (1.5, 1);

			\node at (0.25,0.75) { 1 };
			\node at (0.25,0.25) {\textbf{2}};
			\node at (0.75,0.75) { 3 };
			\node at (0.75,0.25) { 4 };
			\node at (1.25,0.75) { 5 };
			\node at (1.25,0.25) { 6 };
		\end{tikzpicture}
	\end{minipage}
	\begin{minipage}[b]{0.15\textwidth}
		\centering
		\begin{tikzpicture}[line width=1pt]
			\draw[step=0.5] (0,0 ) grid (1.5, 1);

			\node at (0.25,0.75) { 1 };
			\node at (0.25,0.25) {\textbf{2}};
			\node at (0.75,0.75) { 3 };
			\node at (1.25,0.75) { 4 };
			\node at (0.75,0.25) { 5 };
			\node at (1.25,0.25) { 6 };
		\end{tikzpicture}
	\end{minipage}
	\begin{minipage}[b]{0.15\textwidth}
		\centering
		\begin{tikzpicture}[line width=1pt]
			\draw[step=0.5] (0,0 ) grid (1.5, 1);

			\node at (0.25,0.75) { 1 };
			\node at (0.75,0.75) { 2 };
			\node at (0.25,0.25) {\textbf{3}};
			\node at (0.75,0.25) { 4 };
			\node at (1.25,0.75) { 5 };
			\node at (1.25,0.25) { 6 };
		\end{tikzpicture}
	\end{minipage}
	\begin{minipage}[b]{0.15\textwidth}
		\centering
		\begin{tikzpicture}[line width=1pt]
			\draw[step=0.5] (0,0 ) grid (1.5, 1);

			\node at (0.25,0.75) { 1 };
			\node at (0.75,0.75) { 2 };
			\node at (0.25,0.25) {\textbf{3}};
			\node at (1.25,0.75) { 4 };
			\node at (0.75,0.25) { 5 };
			\node at (1.25,0.25) { 6 };
		\end{tikzpicture}
	\end{minipage}
	\begin{minipage}[b]{0.15\textwidth}
		\centering
		\begin{tikzpicture}[line width=1pt]
			\draw[step=0.5] (0,0 ) grid (1.5, 1);

			\node at (0.25,0.75) { 1 };
			\node at (0.75,0.75) { 2 };
			\node at (1.25,0.75) { 3 };
			\node at (0.25,0.25) {\textbf{4}};
			\node at (0.75,0.25) { 5 };
			\node at (1.25,0.25) { 6 };
		\end{tikzpicture}
	\end{minipage}

		\caption{Standard Young tableaux of shape $(2,2)$ and $(3,3)$}
		\label{fig:T3}
\end{figure}

To find the mean of the lower-left entry, it is enough to compute the sum of all lower-left entries. We denote this sum by $\Sigma_n$. Clearly, $\Sigma_1=2$, $\Sigma_2=2+3=5$ and $\Sigma_3=2+2+3+3+4=14$. These values suggest the identity $\Sigma_n=C_{n+1}$, which we prove next.

\begin{proposition}\label{prop:two-row-corner-sum}
For every $n\ge1$,
\[
\Sigma_n
=
\sum_{T\in\SYT(n,n)}T(2,1)
=
C_{n+1}.
\]
\end{proposition}

\begin{figure}[!ht]
\centering
\begin{minipage}[b]{0.3\textwidth}
	\centering
	\begin{tikzpicture}[line width=1pt]
		\draw (0,0) rectangle (2,1);
		\draw (0.5,0) -- (0.5,1);
		\draw (0,0.5) -- (2,0.5);
		\node at (0.25,0.75) {$1$ };
		\node at (0.25,0.25) { $z$};
		\node at (1,0.75) { $\cdots$};
		\node at (1,0.25) { $\cdots$};
	\end{tikzpicture}
\end{minipage}
\begin{minipage}[b]{0.3\textwidth}
	\centering
	\begin{tikzpicture}[line width=1pt]
		\draw (0.5,0) rectangle (2.5,1);
		\draw (1,0) -- (1,1);
		\draw (0.5,0.5) -- (2.5,0.5);
		\draw (0,0) rectangle (0.5,0.5);
		\node at (0.25,0.25) { $w$};
		\node at (1.5,0.75) { $\cdots$};
		\node at (1.5,0.25) { $\cdots$};
	\end{tikzpicture}
\end{minipage}
\caption{Standard Young tableaux of shape $(n,n)$ and of skew shape $(n+1,n+1)/(1)$}
\label{fig:Adding1}
\end{figure}
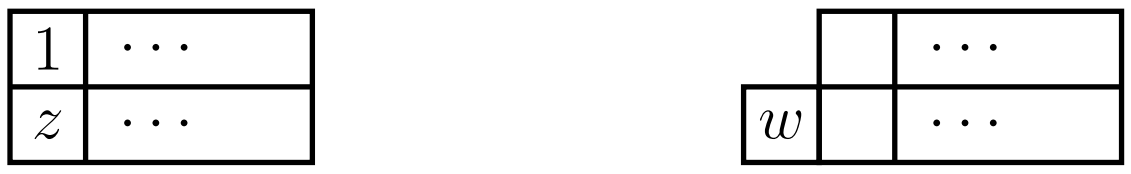

\begin{proof}
Suppose that in a standard Young tableau of shape $(n,n)$, the lower-left entry is $z$, as shown in Figure~\ref{fig:Adding1}. Add one more cell $w$ to the left of $z$, as shown in the same figure. From the figure, the new tableau has skew shape $(n+1,n+1)/(1)$. For each $w\in\{1,\ldots,z\}$, increase every old entry greater than or equal to $w$ by $1$, and fill the new cell with $w$. The new cell has no cell above it, and the entry to its right is $z+1>w$, so the resulting filling is a standard Young tableau of skew shape $(n+1,n+1)/(1)$. Conversely, deleting the new cell and standardizing recovers the original tableau, and the entry $w$ of the deleted cell is at most the entry $z$ to its right. Thus each tableau with lower-left entry $z$ gives exactly $z$ tableaux of the enlarged skew shape, and every tableau of the enlarged skew shape arises exactly once. Moreover, the number of tableaux of skew shape $(n+1,n+1)/(1)$ is the $(n+1)$-th Catalan number: in every standard Young tableau of shape $(n+1,n+1)$ the entry $1$ lies in the cell $(1,1)$, so deleting this cell and decreasing the other entries by $1$ is a bijection onto the tableaux of skew shape $(n+1,n+1)/(1)$. Hence
\[
\Sigma_n=f\bigl((n+1,n+1)/(1)\bigr)=f(n+1,n+1)=C_{n+1}.
\qedhere
\]
\end{proof}

The result above states that the mean of the lower-left entry of a randomly selected standard Young tableau of shape $(n,n)$ is equal to
\[
\frac{C_{n+1}}{C_{n}}=\frac{4n+2}{n+2}.
\]
This gives a combinatorial interpretation of the ratio of two adjacent Catalan numbers. As $n$ goes to infinity, the mean tends to $4$.

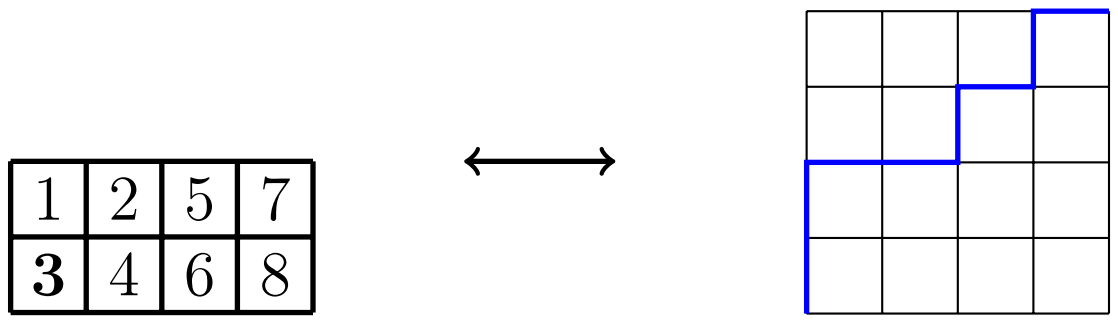
\begin{figure}[!ht]
	\centering
	\begin{minipage}[b]{0.3\textwidth}
		\centering
			\begin{tikzpicture}[line width=1pt]
				\draw[step=0.5] (0,0 ) grid ( 2, 1);
				\node at (0.25,0.75) { 1 };
				\node at (0.75,0.75) { 2 };
				\node at (0.25,0.25) {\textbf{3}};
				\node at (0.75,0.25) { 4 };
				\node at (1.25,0.75) { 5 };
				\node at (1.25,0.25) { 6 };
				\node at (1.75,0.75) { 7 };
				\node at (1.75,0.25) { 8 };

				\draw [<->] (3,1) -- (4,1);
			\end{tikzpicture}
	\end{minipage}
	\begin{minipage}[b]{0.2\textwidth}
		\centering
		\begin{tikzpicture}[line width=1pt]
			\draw[step=0.5,thin] (0,0 ) grid ( 2, 2);

			\draw[blue] (0,0) -- (0,1) -- (1, 1) -- (1,1.5) -- (1.5,1.5) -- (1.5,2) -- (2,2);
		\end{tikzpicture}
	\end{minipage}

	\caption{The bijection between standard Young tableaux of shape $(4,4)$ and Dyck paths of semilength $4$}
	\label{fig:Dyck}
\end{figure}

\begin{remark}\label{rem:dyck}
A \emph{Dyck path} of semilength $n$ is a lattice path from $(0,0)$ to $(n,n)$ with $n$ north steps and $n$ east steps that never goes below the diagonal. There is a bijection between standard Young tableaux of shape $(n,n)$ and Dyck paths of semilength $n$ (see Figure~\ref{fig:Dyck}): if $k$ is in the first row of the tableau, then the $k$-th step of the Dyck path is a north step, and if $k$ is in the second row, then the $k$-th step is an east step. Under this bijection, the lower-left entry $T(2,1)=k$ corresponds to the first east step of the Dyck path occurring at the $k$-th step. Let $B(n,k)$ denote the number of Dyck paths of semilength $n$ whose first east step is the $k$-th step, that is, whose first peak has height $k-1$. The numbers $B(n,k)$ are ballot numbers; see Deutsch~\cite{deutsch1999dyck} for this and other statistics on Dyck paths, and sequence A033184 in the OEIS~\cite{OEIS}. Then Proposition~\ref{prop:two-row-corner-sum} says that
\[
\sum_{k=2}^{n+1}k\,B(n,k)=C_{n+1}.
\]
\end{remark}

The construction in the proof of Proposition~\ref{prop:two-row-corner-sum} can be generalized to standard Young tableaux of arbitrary shape, which gives a skew-shape tableau enumeration formula for the mean of the lower-left entry.

\begin{theorem}\label{thm:general-shape}
Let $\lambda=(\lambda_1,\ldots,\lambda_k)$ be a partition with $k\ge1$ parts, and let $\lambda^{+}=(\lambda_1+1,\ldots,\lambda_k+1)$. Then
\begin{equation}\label{eq:general-shape}
\sum_{T\in\SYT(\lambda)}T(k,1)=f(\lambda^{+}/1^{k-1}).
\end{equation}
In particular, for a uniformly random standard Young tableau of shape $\lambda$, the mean of the lower-left entry is $f(\lambda^{+}/1^{k-1})/f(\lambda)$.
\end{theorem}

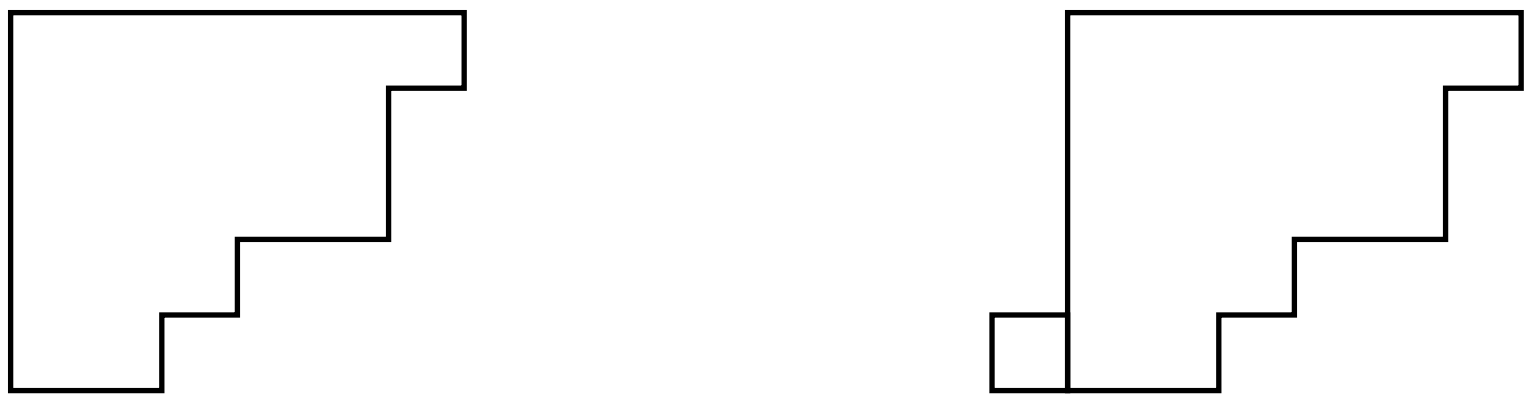
\begin{figure}[!ht]
	\centering
	\begin{minipage}[b]{0.4\textwidth}
		\centering
		\begin{tikzpicture}[line width=1pt]
			\draw (0,0) -- (0,2.5) -- (3,2.5) -- (3,2) -- (2.5,2) -- (2.5,1) -- (1.5,1) -- (1.5,0.5) -- (1,0.5) -- (1,0) -- cycle;
		\end{tikzpicture}
	\end{minipage}
	\begin{minipage}[b]{0.4\textwidth}
		\centering
		\begin{tikzpicture}[line width=1pt]
			\draw (0,0) -- (0,2.5) -- (3,2.5) -- (3,2) -- (2.5,2) -- (2.5,1) -- (1.5,1) -- (1.5,0.5) -- (1,0.5) -- (1,0) -- cycle;
			\draw (0,0) -- (-0.5,0) -- (-0.5,0.5) -- (0, 0.5) -- cycle;
		\end{tikzpicture}
	\end{minipage}

	\caption{A shape $\lambda$ and the skew shape $\lambda^{+}/1^{k-1}$}
	\label{fig:Adding4}
\end{figure}

\begin{proof}
Add a new cell immediately to the left of the lower-left cell of $\lambda$. The resulting skew shape is $\lambda^{+}/1^{k-1}$, as shown in Figure~\ref{fig:Adding4}. Exactly as in the proof of Proposition~\ref{prop:two-row-corner-sum}, a tableau $T\in\SYT(\lambda)$ with lower-left entry $z$ gives exactly $z$ standard Young tableaux of shape $\lambda^{+}/1^{k-1}$, one for each $w\in\{1,\ldots,z\}$, and every tableau of shape $\lambda^{+}/1^{k-1}$ arises exactly once. Therefore the sum of all lower-left entries over $\SYT(\lambda)$ is $f(\lambda^{+}/1^{k-1})$. Dividing by the total number $f(\lambda)$ of standard Young tableaux of shape $\lambda$ gives the mean.
\end{proof}

We first apply Theorem~\ref{thm:general-shape} to two-row shapes; the resulting formula is the case $\lambda=(m,n)$ of \eqref{eq:RN}.

\begin{corollary}\label{cor:two-row-mean}
Let $m\ge n\ge1$. Then
\[
\sum_{T\in\SYT(m,n)}T(2,1)=f(m+1,n+1),
\]
and for $T$ chosen uniformly from $\SYT(m,n)$,
\[
\Ex[T(2,1)]
=
\frac{(m+n+2)(m+n+1)}{(m+2)(n+1)}.
\]
\end{corollary}

\begin{proof}
By Theorem~\ref{thm:general-shape}, the sum equals $f((m+1,n+1)/(1))$, which equals $f(m+1,n+1)$ as in the proof of Proposition~\ref{prop:two-row-corner-sum}. It then follows from \eqref{eq:two-row-f} that
\begin{align*}
\Ex[T(2,1)]
&=\frac{f(m+1,n+1)}{f(m,n)}
=\frac{m-n+1}{m+2}\binom{m+n+2}{n+1}\Big/\frac{m-n+1}{m+1}\binom{m+n}{n}\\
&=\frac{m+1}{m+2}\cdot\frac{(m+n+2)(m+n+1)}{(n+1)(m+1)}
=\frac{(m+n+2)(m+n+1)}{(m+2)(n+1)}.
\qedhere
\end{align*}
\end{proof}

For rectangular shapes, we also need the following lemma on rotating a tableau by $180^\circ$.

\begin{lemma}\label{lem:rotate-complement}
Let $\lambda=(\lambda_1,\lambda_2,\ldots,\lambda_k)\vdash n$. Define
$\mu=(\lambda_1)^k$ and
\[
\nu=(\lambda_1-\lambda_k,\lambda_1-\lambda_{k-1},\ldots,\lambda_1-\lambda_1).
\]
Then
\[
f(\lambda)=f(\mu/\nu).
\]
\end{lemma}

\begin{proof}
Rotate a tableau $T\in\SYT(\lambda)$ by $180^\circ$. If the entry $T(i,j)$ is in the cell $(i,j)$, then after rotation it is placed in the cell
\[
i'=k+1-i,\qquad j'=\lambda_1+1-j,
\]
and its value is changed to $n+1-T(i,j)$. This operation converts the increasing row and column conditions of $T$ into the increasing row and column conditions of a standard Young tableau of skew shape $\mu/\nu$. The construction is reversible, hence it gives a bijection between $\SYT(\lambda)$ and $\SYT(\mu/\nu)$.
\end{proof}

\begin{corollary}\label{cor:rect}
For $n,k\ge1$,
\begin{equation}\label{eq:rect-sum}
\sum_{T\in\SYT(n^k)}T(k,1)
=
f(n+1,n^{k-1}),
\end{equation}
and for a uniformly random standard Young tableau of shape $n^k$, the mean of the lower-left entry is
\begin{equation}\label{eq:rect-mean}
\frac{k(nk+1)}{n+k}=\frac{k^2n+k}{n+k}.
\end{equation}
For fixed $k$, the mean tends to $k^2$ as $n\to\infty$.
\end{corollary}

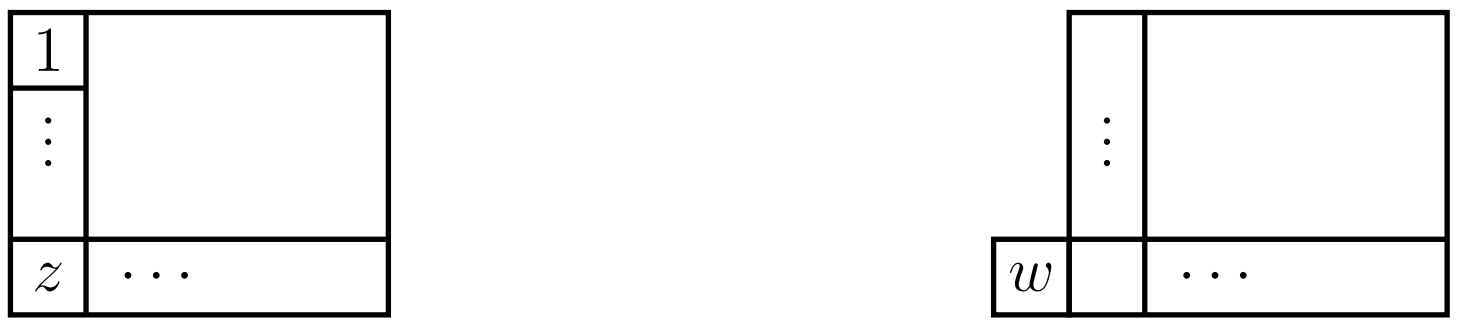
\begin{figure}[!ht]
	\centering
	\begin{minipage}[b]{0.4\textwidth}
		\centering
		\begin{tikzpicture}[line width=1pt]
			\draw (0,0) rectangle (2.5,2);
			\draw (0.5,0) -- (0.5,2);
			\draw (0,0.5) -- (2.5,0.5);
			\draw (0,1.5) -- (0.5,1.5);
			\node at (0.25,0.25) { $z$};
			\node at (0.25,1.75) { $1$};
			\node at (0.25,1.25) { $\vdots$};
			\node at (1,0.25) { $\cdots$};
		\end{tikzpicture}
	\end{minipage}
	\begin{minipage}[b]{0.4\textwidth}
		\centering
		\begin{tikzpicture}[line width=1pt]
			\draw (1,0) rectangle (3.5,2);
			\draw (1.5,0) -- (1.5,2);
			\draw (1,0.5) -- (3.5,0.5);
			\draw (0.5,0) rectangle (1,0.5);
			\node at (0.75,0.25) { $w$};
			\node at (1.25,1.25) { $\vdots$};
			\node at (2,0.25) { $\cdots$};
		\end{tikzpicture}
	\end{minipage}

	\caption{Standard Young tableaux of shape $n^k$ and of skew shape $(n+1)^k/1^{k-1}$}
	\label{fig:Adding2}
\end{figure}

\begin{proof}
By Theorem~\ref{thm:general-shape},
\[
\sum_{T\in\SYT(n^k)}T(k,1) = f\bigl((n+1)^k/1^{k-1}\bigr);
\]
the skew shape $(n+1)^k/1^{k-1}$ is shown in Figure~\ref{fig:Adding2}. Rotating this skew shape through $180^\circ$ and replacing each entry $t$ by $nk+2-t$ gives a bijection with standard Young tableaux of shape $(n+1,n^{k-1})$; this is Lemma~\ref{lem:rotate-complement} applied to $\lambda=(n+1,n^{k-1})$, with $\mu=(n+1)^k$ and $\nu=1^{k-1}$. Hence
\[
f\bigl((n+1)^k/1^{k-1}\bigr) = f(n+1,n^{k-1}),
\]
which is \eqref{eq:rect-sum}. The shapes $n^k$ and $(n+1,n^{k-1})$ have the same hook lengths outside the first row. In the first row, the hook lengths of $n^k$ are $n+k-1,n+k-2,\ldots,k$, while those of $(n+1,n^{k-1})$ are $n+k,n+k-1,\ldots,k+1$ and $1$. Hence, by the hook length formula \eqref{eq:hook},
\begin{align*}
\frac{f(n+1,n^{k-1})}{f(n^k)}
&=\frac{(nk+1)!}{(nk)!}\cdot\frac{(n+k-1)(n+k-2)\cdots(k+1)k}{(n+k)(n+k-1)\cdots(k+1)\cdot1}\\
&=\frac{(nk+1)k}{n+k}
=\frac{k^2n+k}{n+k},
\end{align*}
which is \eqref{eq:rect-mean}; letting $n$ tend to infinity gives the limit $k^2$.
\end{proof}

For $k=3$, the right side of \eqref{eq:rect-sum} is the number of standard Young tableaux of shape $(n+1,n,n)$, which is $3,21,210,2574,\ldots$ for $n=1,2,3,4,\ldots$; see sequence A123691 in the OEIS~\cite{OEIS}. Thus this sequence also gives the sum of the lower-left entries over $\SYT(n,n,n)$.

\subsection{Distributions in two-row shapes}\label{subsec:two-row-distributions}

We next determine the distribution of the entry in an arbitrary cell of the second row of a two-row tableau. Ekhad and Zeilberger~\cite{ekhad2024experimenting} compute the distribution of the entry in a given cell of a random tableau of a given shape by decomposing the tableau according to the subshape occupied by the entries up to a prescribed value, and summing over all possible subshapes. For the cell $(2,k)$ of a two-row shape this subshape is unique, which gives the product formula \eqref{eq:dist-f} below, and we then evaluate the product explicitly by the determinantal formula \eqref{eq:aitken}.

\begin{theorem}\label{thm:two-row-dist}
Let $m\ge n\ge1$ and $1\le k\le n$, and let $T$ be chosen uniformly from $\SYT(m,n)$. Then $2k\le T(2,k)\le m+k$, and for $2k\le a\le m+k$,
\begin{align}
\Pr(T(2,k)=a)
&=
\frac{f(a-k,k-1)\,f\bigl((m,n)/(a-k,k)\bigr)}{f(m,n)}
\label{eq:dist-f}\\
&=
\frac{
\left[
\dbinom{a-1}{k-1}-\dbinom{a-1}{k-2}
\right]
\left[
\dbinom{m+n-a}{n-k}
-\dbinom{m+n-a}{n-a+k-1}
\right]
}{
\dbinom{m+n}{n}-\dbinom{m+n}{n-1}
}.
\label{eq:dist-explicit}
\end{align}
In particular, when $a=2k$, the entry $T(2,k)$ takes its smallest possible value and
\begin{equation}\label{eq:dist-min}
\Pr(T(2,k)=2k)
=
C_k\cdot
\frac{m+1}{m-k+1}\cdot
\frac{\dbinom{m+n-2k}{n-k}}{\dbinom{m+n}{n}},
\end{equation}
where $C_k=\binom{2k}{k}/(k+1)$ is the $k$-th Catalan number.
\end{theorem}

\begin{proof}
Suppose that $T(2,k)=a$. As illustrated in Figure~\ref{fig:split}, the entries less than $a$ fill the blue region, which consists of the cells $(2,1),\ldots,(2,k-1)$ of the second row and the cells $(1,1),\ldots,(1,a-k)$ of the first row. Since $T(1,k)<a$, we have $a-k\ge k$, and since $a-k\le m$, we have $a\le m+k$. Thus $2k\le a\le m+k$. The total number of fillings of the blue region is $f(a-k,k-1)$. The cells in the orange region are filled with entries greater than $a$, and the number of admissible fillings is the number of standard Young tableaux of the skew shape $(m,n)/(a-k,k)$. Multiplying the two counts and dividing by $f(m,n)$ gives \eqref{eq:dist-f}.

\begin{figure}[!ht]
	\centering
	\begin{minipage}[b]{0.45\textwidth}
		\centering
		\begin{tikzpicture}[line width=1pt]
			\draw (0,0) rectangle (3,0.5);
			\draw (0,0.5) rectangle (4.5,1);
			\draw[blue] (0,0) -- (0,1) -- (2.5,1) -- (2.5,0.5) -- (1,0.5) -- (1,0) -- cycle;
			\draw[orange] (1.5,0) -- (1.5,0.5) -- (2.5,0.5) -- (2.5,1) -- (4.5,1) -- (4.5,0.5) -- (3,0.5) -- (3,0) -- cycle;
			\node at (1.25,0.25) { $a$ };
			\node at (0.5,0.75) { $\cdots$ };
			\node at (0.5,0.25) { $\cdots$ };
			\node at (3,0.75) { $\cdots$ };
			\node at (2,0.25) { $\cdots$ };
			\node at (2.25,-0.5) {(a) The case $n>a-k$};
		\end{tikzpicture}
	\end{minipage}
	\begin{minipage}[b]{0.45\textwidth}
		\centering
		\begin{tikzpicture}[line width=1pt]
			\draw (0,0) rectangle (3,0.5);
			\draw (0,0.5) rectangle (4.5,1);
			\draw[blue] (0,0) -- (0,1) -- (3.5,1) -- (3.5,0.5) -- (1,0.5) -- (1,0) -- cycle;
			\draw[orange] (1.5,0) rectangle (3,0.5);
			\draw[orange] (3.5,0.5) rectangle (4.5,1);
			\node at (1.25,0.25) { $a$};
			\node at (0.5,0.75) { $\cdots$};
			\node at (0.5,0.25) { $\cdots$};
			\node at (4,0.75) { $\cdots$};
			\node at (2,0.25) { $\cdots$};
			\node at (2.25,-0.5) {(b) The case $n\leq a-k$};
		\end{tikzpicture}
	\end{minipage}
	\caption{A standard Young tableau of shape $(m,n)$ with $T(2,k)=a$: the entries smaller than $a$ fill the blue region, the entries greater than $a$ fill the orange region}
	\label{fig:split}
\end{figure}
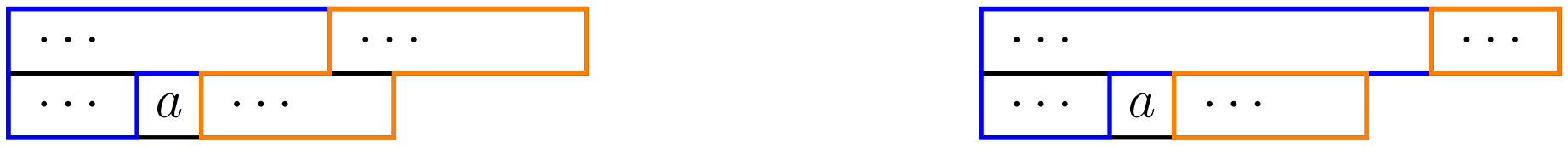

By \eqref{eq:aitken} with $\lambda=(m,n)$ and $\mu=(a-k,k)$,
\begin{align*}
f\bigl((m,n)/(a-k,k)\bigr)
&=(m+n-a)!
\begin{vmatrix}
\dfrac{1}{(m-a+k)!} & \dfrac{1}{(m-k+1)!}\\[1.2ex]
\dfrac{1}{(n-a+k-1)!} & \dfrac{1}{(n-k)!}
\end{vmatrix}\\
&=\binom{m+n-a}{n-k}
-\binom{m+n-a}{n-a+k-1}.
\end{align*}
By \eqref{eq:two-row-f},
\[
f(a-k,k-1)
=\binom{a-1}{k-1}-\binom{a-1}{k-2}
=\frac{a-2k+2}{a-k+1}\binom{a-1}{k-1}
\]
and
\[
f(m,n)
=\binom{m+n}{n}-\binom{m+n}{n-1}
=\frac{m-n+1}{m+1}\binom{m+n}{n}.
\]
Substituting these three expressions into \eqref{eq:dist-f} gives \eqref{eq:dist-explicit}.

When $a=2k$, we have
\[
\binom{2k-1}{k-1}-\binom{2k-1}{k-2}=C_k
\]
and
\[
\binom{m+n-2k}{n-k}-\binom{m+n-2k}{n-k-1}
=\frac{m-n+1}{m-k+1}\binom{m+n-2k}{n-k},
\]
which yields the stated special case \eqref{eq:dist-min}.
\end{proof}

\subsection{Sums of second-row entries}\label{subsec:second-row-sums}

For the shape $(n,n)$, the insertion idea, phrased in terms of Dyck paths, also gives the sum of the entries in every cell of the second row. Under the bijection described in Remark~\ref{rem:dyck}, $T(2,k)$ is the position of the $k$-th east step of the corresponding Dyck path. The \emph{height} of a Dyck path after a step is the number of north steps minus the number of east steps taken so far, and the path \emph{returns to the diagonal} at a step after which its height is $0$.

\begin{theorem}\label{thm:second-row-sums}
Let $1\le k\le n$, and put $T(2,0)=0$. Then
\begin{equation}\label{eq:increment}
\sum_{T\in\SYT(n,n)}\bigl(T(2,k)-T(2,k-1)\bigr)=\sum_{j=k-1}^{n}C_jC_{n-j}.
\end{equation}
Consequently,
\begin{equation}\label{eq:second-row-sum}
\sum_{T\in\SYT(n,n)}T(2,k)
=\sum_{j=0}^{n}\min(j+1,k)\,C_jC_{n-j}
=kC_{n+1}-\sum_{i=1}^{k-1}(k-i)\,C_{i-1}C_{n+1-i}.
\end{equation}
\end{theorem}

For $k=1$, \eqref{eq:second-row-sum} is Proposition~\ref{prop:two-row-corner-sum}. For $k=2$ and $n\ge2$, it gives
\[
\sum_{T\in\SYT(n,n)}T(2,2)=2C_{n+1}-C_n.
\]
The numbers $2C_{n+1}-C_n$ form sequence A000782 in the OEIS~\cite{OEIS}.

\begin{proof}
Let $P$ be the Dyck path corresponding to $T$, and let $E_i$ be the position of its $i$-th east step, with $E_0=0$. Then $T(2,k)-T(2,k-1)=E_k-E_{k-1}$ is the number of integers $w$ satisfying $E_{k-1}<w\le E_k$. These indices correspond to the north steps preceding the $k$-th east step and following the previous east step, if any, together with the $k$-th east step itself. Thus the left side of \eqref{eq:increment} counts the pairs $(P,w)$, where $P$ is a Dyck path of semilength $n$ and $w$ is an integer satisfying $E_{k-1}<w\le E_k$.

Given such a pair, write $P=P_1P_2\cdots P_{2n}$ and let
\[
Q=N\,P_1\cdots P_{w-1}\,E\,P_w\cdots P_{2n},
\]
that is, insert a north step at the beginning of $P$ and an east step immediately before the step $P_w$. The steps $P_1,\ldots,P_{w-1}$ contain exactly $k-1$ east steps, so the inserted east step is the $k$-th east step of $Q$. Before the inserted east step, $Q$ runs one unit above $P$ and hence stays at height at least $1$; after it, $Q$ runs at the same height as $P$. Therefore $Q$ is a Dyck path of semilength $n+1$ that does not return to the diagonal before its $k$-th east step.

Conversely, let $Q$ be a Dyck path of semilength $n+1$ that does not return to the diagonal before its $k$-th east step. Then $Q$ starts with a north step, and before its $k$-th east step it stays at height at least $1$. Deleting the first step and the $k$-th east step of $Q$ gives a Dyck path $P$ of semilength $n$. Since $k\le n$, some step of $Q$ follows the deleted east step; its position $w$ in $P$ satisfies $E_{k-1}<w\le E_k$, because the steps before it contain exactly $k-1$ east steps. This reverses the construction above.

Finally, let $Q$ be a Dyck path of semilength $n+1$, and let its first return to the diagonal occur at its $(j+1)$-th east step, where $0\le j\le n$. Then $Q=N\,P'\,E\,P''$ with Dyck paths $P'$ and $P''$ of semilengths $j$ and $n-j$, and $Q$ does not return to the diagonal before its $k$-th east step exactly when $j+1\ge k$. Hence the number of such $Q$ is $\sum_{j=k-1}^{n}C_jC_{n-j}$, which proves \eqref{eq:increment}.

Applying \eqref{eq:increment} with $k$ replaced by $i$ and summing over $i=1,\ldots,k$, we see that, for each fixed $j$, the term $C_jC_{n-j}$ occurs once for each integer $i$ satisfying $1\le i\le k$ and $i-1\le j$, hence $\min(j+1,k)$ times. This gives the first expression in \eqref{eq:second-row-sum}. Since $\sum_{j=0}^{n}C_jC_{n-j}=C_{n+1}$,
\[
\sum_{j=0}^{n}\min(j+1,k)\,C_jC_{n-j}
=kC_{n+1}-\sum_{j=0}^{k-2}(k-1-j)\,C_jC_{n-j},
\]
and the substitution $i=j+1$ gives the second expression.
\end{proof}

Dividing \eqref{eq:second-row-sum} by $C_n$ and letting $n\to\infty$, we recover a limit formula found by Stanley; see Chan, Pak and Panova~\cite{chan2021sorting}.

\begin{corollary}\label{cor:stanley-limit}
For fixed $k\ge1$,
\[
\lim_{n\to\infty}\frac{1}{C_n}\sum_{T\in\SYT(n,n)}T(2,k)=2k+\frac{k}{4^{k-1}}\binom{2k}{k}.
\]
\end{corollary}

\begin{proof}
By \eqref{eq:second-row-sum},
\[
\frac{1}{C_n}\sum_{T\in\SYT(n,n)}T(2,k)
=k\,\frac{C_{n+1}}{C_n}-\sum_{i=1}^{k-1}(k-i)\,C_{i-1}\,\frac{C_{n+1-i}}{C_n}.
\]
Since $C_{m+1}/C_m=(4m+2)/(m+2)$ tends to $4$, we have $C_{n+1-i}/C_n\to4^{1-i}$ for each fixed $i$. Hence the limit, which we denote by $L_k$, exists and
\[
L_k=4k-\sum_{i=1}^{k-1}(k-i)\,C_{i-1}\,4^{1-i}.
\]
We compute the generating function of the numbers $L_k$. By the generating function
\[
\sum_{j\ge0}C_jx^j=\frac{1-\sqrt{1-4x}}{2x}
\]
of the Catalan numbers, with $x$ replaced by $x/4$,
\[
\sum_{i\ge1}C_{i-1}4^{1-i}x^i=x\sum_{j\ge0}C_j\left(\frac{x}{4}\right)^{j}=2\bigl(1-\sqrt{1-x}\bigr).
\]
Since $\sum_{m\ge1}m\,x^m=x/(1-x)^2$, it follows that
\[
\sum_{k\ge1}L_k\,x^k
=\frac{4x}{(1-x)^2}-\frac{2x\bigl(1-\sqrt{1-x}\bigr)}{(1-x)^2}
=\frac{2x}{(1-x)^2}+\frac{2x}{(1-x)^{3/2}}.
\]
The first term is $\sum_{k\ge1}2k\,x^k$. For the second term, we differentiate the binomial series of $(1-x)^{-1/2}$:
\begin{align*}
\frac{2x}{(1-x)^{3/2}}
&=4x\,\frac{d}{dx}\,\frac{1}{\sqrt{1-x}}
=4x\,\frac{d}{dx}\sum_{m\ge0}\binom{2m}{m}\frac{x^m}{4^m}\\
&=\sum_{k\ge1}\frac{4k}{4^{k}}\binom{2k}{k}x^k
=\sum_{k\ge1}\frac{k}{4^{k-1}}\binom{2k}{k}x^k.
\end{align*}
Comparing the coefficients of $x^k$ gives the stated formula.
\end{proof}

\begin{remark}
The same bijection applies to a two-row shape $(m,n)$ with $m\ge n\ge k\ge1$. The tableaux of shape $(m,n)$ correspond to lattice paths with $m$ north steps and $n$ east steps that never go below the diagonal. A path for the shape $(m+1,n+1)$ either returns to the diagonal, in which case it has the first return decomposition $N\,P'\,E\,P''$ with $P'$ a Dyck path of semilength $j$ and $P''$ a path for the shape $(m-j,n-j)$, or never returns to the diagonal, in which case $m>n$ and it has the form $N\,P'$ with $P'$ a path for the shape $(m,n+1)$. The argument in the proof of Theorem~\ref{thm:second-row-sums} therefore gives
\[
\sum_{T\in\SYT(m,n)}\bigl(T(2,k)-T(2,k-1)\bigr)=f(m,n+1)+\sum_{j=k-1}^{n}C_j\,f(m-j,n-j),
\]
where $f(m,n+1)$ is given by \eqref{eq:two-row-f} and is $0$ when $m=n$. For $k=1$ the right side is $f(m+1,n+1)$ by the first return decomposition, and the formula is the first identity of Corollary~\ref{cor:two-row-mean}.
\end{remark}

\subsection{Pattern-avoiding alternating permutations}\label{subsec:pattern-avoidance}

We conclude this section by applying the results for the shape $(n,n)$ to pattern-avoiding alternating permutations. For $T\in\SYT(n,n)$, let $w(T)$ be the word obtained by reading the columns of $T$ from left to right, each column from bottom to top:
\[
w(T)=T(2,1)\,T(1,1)\,T(2,2)\,T(1,2)\cdots T(2,n)\,T(1,n).
\]
For example, the tableau in Figure~\ref{fig:Dyck} gives $w(T)=31426587$.

\begin{lemma}\label{lem:reading-word}
The map $T\mapsto w(T)$ is a bijection from $\SYT(n,n)$ onto $\cA'_{2n}(321)$. Moreover, $w(T)_{2k-1}=T(2,k)$ and $w(T)_{2k}=T(1,k)$ for $1\le k\le n$.
\end{lemma}

\begin{proof}
Let $w=w(T)$. Since the columns of $T$ increase, $w_{2k-1}>w_{2k}$, and since $T(1,k)<T(1,k+1)<T(2,k+1)$, we have $w_{2k}<w_{2k+1}$. Thus $w$ is down-up. The entries of $w$ in odd positions form the second row of $T$ and those in even positions form the first row, so $w$ is a union of two increasing subsequences and avoids $321$.

Conversely, let $w\in\cA'_{2n}(321)$. If $w_{2k-1}>w_{2k+1}$ for some $k<n$, then $w_{2k-1}>w_{2k+1}>w_{2k+2}$ is an occurrence of $321$. If $w_{2k}>w_{2k+2}$ for some $k<n$, then $w_{2k-1}>w_{2k}>w_{2k+2}$ is an occurrence of $321$. Hence the entries of $w$ in odd positions increase, and so do the entries in even positions. Putting $T(2,k)=w_{2k-1}$ and $T(1,k)=w_{2k}$ therefore gives a standard Young tableau $T$ of shape $(n,n)$ with $w(T)=w$.
\end{proof}

Up to reversing the permutation and transposing the tableau, the map $T\mapsto w(T)$ is the alternating case of a bijection of Lewis~\cite{lewis2011pattern} onto standard Young tableaux of rectangular shape. In particular, $|\cA'_{2n}(321)|=C_n$. By complementation, $|\cA_{2n}(123)|=C_n$ as well, a result that Stanley attributes to Deutsch and Reifegerste; see Lewis~\cite{lewis2011pattern}.

\begin{theorem}\label{thm:pattern}
Let $w$ be chosen uniformly from $\cA'_{2n}(321)$ and let $1\le k\le n$. Then
\[
\Ex[w_{2k-1}]
=\frac{1}{C_n}\left(kC_{n+1}-\sum_{i=1}^{k-1}(k-i)\,C_{i-1}C_{n+1-i}\right)
\]
and $\Ex[w_{2k}]=2n+1-\Ex[w_{2n+1-2k}]$. In particular, $\Ex[w_1]=C_{n+1}/C_n$. Moreover, $2k\le w_{2k-1}\le n+k$, and for $2k\le a\le n+k$,
\[
\Pr(w_{2k-1}=a)
=\frac{1}{C_n}\left[\binom{a-1}{k-1}-\binom{a-1}{k-2}\right]\left[\binom{2n-a}{n-k}-\binom{2n-a}{n-a+k-1}\right].
\]
Finally, for fixed $k$,
\[
\lim_{n\to\infty}\Ex[w_{2k-1}]=2k+\frac{k}{4^{k-1}}\binom{2k}{k}.
\]
\end{theorem}

\begin{proof}
By Lemma~\ref{lem:reading-word}, $w_{2k-1}$ has the same distribution as $T(2,k)$ for $T$ chosen uniformly from $\SYT(n,n)$. The statements on $w_{2k-1}$ therefore follow from Theorem~\ref{thm:second-row-sums}, Theorem~\ref{thm:two-row-dist} with $m=n$, and Corollary~\ref{cor:stanley-limit}. Reversing a permutation and replacing each entry $v$ by $2n+1-v$ maps $\cA'_{2n}(321)$ to itself and sends the entry in position $2k$ to $2n+1$ minus the entry in position $2n+1-2k$. This gives $\Ex[w_{2k}]=2n+1-\Ex[w_{2n+1-2k}]$.
\end{proof}

Since $2n+1-2k=2(n+1-k)-1$ is odd, the theorem gives the mean of every entry of $w$. For $k=1$, the number of $w\in\cA'_{2n}(321)$ with $w_1=a$ is the ballot number $B(n,a)$ of Remark~\ref{rem:dyck}; compare the result on the first entry of $123$- and $132$-avoiding permutations of Miner and Pak~\cite{miner2014shape} mentioned in the introduction. By complementation, the theorem also describes $123$-avoiding up-down permutations: if $\sigma$ is chosen uniformly from $\cA_{2n}(123)$, then $\sigma_i$ has the same distribution as $2n+1-w_i$. In particular, $\Ex[\sigma_1]=2n+1-C_{n+1}/C_n$, while Theorem~\ref{thm:first-entry} gives $\Ex[\sigma_1]=2n+1-E_{2n+1}/E_{2n}$ for $\sigma$ chosen uniformly from all of $\cA_{2n}$.

\section{Concluding remarks}
This paper has studied the first entry and the positions of fixed values in alternating permutations by using linear extensions and generating functions. For standard Young tableaux, an insertion construction gives a skew-tableau interpretation of the sum of the lower-left entries, with explicit formulas for rectangular shapes. For two-row shapes, we also gave the distribution of the entry in each cell of the second row and, for the shape $(n,n)$, a closed formula for the sum of the entries in each cell of the second row; through a bijection of Lewis, these results describe the entries of a random $321$-avoiding down-up permutation.

Three natural questions arise. First, for shapes with three or more rows, such as $(n,n,n)$, we do not know formulas for the sums of the entries $T(i,k)$ with $i,k\ge2$. Second, Lewis~\cite{lewis2011pattern} showed that the up-down permutations of length $2n$ avoiding $1234$ are equinumerous with the standard Young tableaux of shape $(n,n,n)$. It would be interesting to find a bijection under which the lower-left entry of the tableau, whose sum is given by Corollary~\ref{cor:rect} with $k=3$, corresponds to a natural statistic of the permutation. Third, it would be interesting to extend the linear extension method to related families of posets.

\par\medskip\noindent\textbf{AI declaration.} During the preparation of this work the authors used GPT-5.6 Sol to improve the language of the manuscript. After using this tool, the authors reviewed and edited the content as needed and take full responsibility for the content of the publication.

\end{document}